\documentclass[12pt]{amsart}

\newtheorem{theorem}{Theorem}[section] 
\newtheorem{corollary}[theorem]{Corollary} 
\newtheorem{conjecture}[theorem]{Conjecture} 
\newtheorem{lemma}[theorem]{Lemma} 
\newtheorem{proposition}[theorem]{Proposition}

\theoremstyle{definition}

\theoremstyle{remark}
\newtheorem{remark}[equation]{Remark}

  {\end{list}}

\usepackage{amssymb} 
\usepackage{verbatim}

\def\ord{{\mathop{\rm ord}}}
\def\Gal{{\mathop{\rm Gal}}}

\def\supp{{\mathop{\rm supp}}}

\def\deg{{\mathop{\rm deg}}}

\def\Spec{{\mathop{\rm Spec}}}

\def\Im{{\mathop{\rm Im}}}

\def\tors{{\mathop{\rm tors}}}
\def\cl{{\mathop{\rm cl}}}
\def\larg{{\mathop{\rm larg}}}
\def\GL{{\mathop{\rm GL}}}
\def\area{{\mathop{\rm area}}}

\def\<{\langle }
\def\>{\rangle }

\def\({(\!(}
\def\){)\!)}

\def\BF1{{\mathbf 1}}

\def\hhat{{\widehat{h}}}

\def\fp{{\mathfrak{p}}}

\def\tlambda{{\tilde{\lambda}}}

\def\AA{{\mathbb A}}
\def\CC{{\mathbb C}}
\def\QQ{{\mathbb Q}}

\def\RR{{\mathbb R}}
\def\ZZ{{\mathbb Z}}
\def\PP{{\mathbb P}}

\def\cA{{\mathcal A}}
\def\cO{{\mathcal O}}

\def\cE{{\mathcal E}}

\def\cF{{\mathcal F}}

\def\cS{{\mathcal S}}

\def\bcO{{\overline{\cO}}}

\def\kbar{{\overline{k}}}

\def\kbar{{\overline{k}}}

\def\Lbar{{\overline{L}}}

\def\Zbar{{\overline{\mathbb Z}}}

\def\Qbar{{\overline{\QQ}}}

\def\ndiv{{\!\not\ \mid}}

\DeclareMathSymbol{\varnothing} {\mathord}{AMSb}{"3F} 
 
\theoremstyle{definition} 
\theoremstyle{remark} 

\begin{document}
\title{A finiteness property of torsion points}

\author{Matthew Baker}
\address{Matthew Baker\\ 
Department of Mathematics\\
Georgia Institute of Technology\\
Atlanta, Georgia 30332\\
USA}
\email{mbaker@@math.gatech.edu}

\author{Su-Ion Ih}
\address{Su-Ion Ih\\ 
Department of Mathematics \\
University of Colorado at Boulder \\
Campus Box 395 \\
Boulder, CO 80309-0395 \\
USA}
\email{Su-ion.Ih@@colorado.edu}

\author{Robert Rumely}
\address{Robert Rumely\\ 
Department of Mathematics\\
University of Georgia\\
Athens, Georgia 30602\\
USA}
\email{rr@@math.uga.edu}

\date{September 19, 2005}
\thanks{Work supported in part by NSF grant DMS-0300784.}

\begin{abstract}
Let $k$ be a number field, let $E/k$ be an elliptic curve, and
let $S$ be a finite set of places of $k$ containing the archimedean places. 
We prove that if  $\alpha \in E(\kbar)$ is nontorsion,
then there are only finitely many torsion points $\xi \in E(\kbar)_{\tors}$ 
which are $S$-integral with respect to $\alpha$. 
We also prove an analogue of this for the multiplicative group,
and formulate conjectural generalizations for abelian
varieties and dynamical systems.
\end{abstract}

\maketitle

\section{Introduction}

Let $k$ be a number field, with ring of integers $\cO_k$ and algebraic closure 
$\kbar$, and let $E/k$ be an elliptic curve.  
Let $\cE/\Spec(\cO_k)$ be a model of $E$,
and let $S$ be a finite set of places of $k$ containing the archimedean places. 
In this paper we will prove:

\begin{theorem} \label{MainThm} 
If  $\alpha \in E(\kbar)$ is nontorsion
$($that is, has canonical height $\hhat(\alpha) > 0$$)$,
then there are only finitely many torsion points $\xi \in E(\kbar)_{\tors}$ 
which are $S$-integral with respect to $\alpha$. 
\end{theorem} 
 
\noindent{By} $S$-integrality we mean that the Zariski closures of 
$\xi$ and $\alpha$ in the model $\cE/\Spec(\cO_{k})$ 
do not meet outside fibres above $S$.   Since any two models 
are isomorphic outside a finite set of places, 
the finiteness property is independent of the set $S$ and the model $\cE$.

We will also prove an analogue of Theorem~\ref{MainThm} 
for the multiplicative group
(Theorem~\ref{Thm1} below). Theorems \ref{MainThm} and \ref{Thm1} 
are analogues for non-compact varieties, 
where it is most natural to look at integral points, 
of the Manin-Mumford conjecture (first proved by Raynaud \cite{Ray}) .  

\vskip .1 in
The ingredients of the proof of Theorem~\ref{MainThm} are a strong form of equidistribution 
for torsion points at all places $v$, properties of local height functions, 
and David/Hirata-Kohno's theorem on linear forms in elliptic logarithms.  
In outline, the proof is as follows.
By base change, one reduces to the case where $\alpha \in E(k)$.  
Given a place $v$ of $k$, 
let $\kbar_v$ be the algebraic closure of the completion $k_v$, 
and let $\lambda_v : E(\kbar_v) \rightarrow \RR$ be an appropriately normalized
N{\'e}ron-Tate canonical local height.  
On the one hand, elementary properties of heights show that
for any torsion point $\xi_n$, one has 
\begin{equation} \label{Intro1} 
\hhat(\alpha) \ = \ \frac{1}{[k(\xi_n):k]} \ \sum_v 
\sum_{\sigma : k(\xi_n)/k \hookrightarrow \kbar_v} \lambda_v(\alpha-\sigma(\xi_n)) \ . 
\end{equation}  
On the other hand, if $\{\xi_n\}$
is a sequence of distinct torsion points which are $S$-integral with respect
to $\alpha$, then for each $v$ 
\begin{equation} \label{Intro2} 
\lim_{n \rightarrow \infty} \frac{1}{[k(\xi_n):k]} \sum_{\sigma : k(\xi_n)/k \hookrightarrow \kbar_v} \lambda_v(\alpha-\sigma(\xi_n)) \ = \ 0 \ .
\end{equation} 
By the integrality hypothesis, the outer sum in (\ref{Intro1}) 
can be restricted to $v \in S$, allowing the limit and the sum to be interchanged.
This gives $\hhat(\alpha) = 0$, contradicting the assumption that $\alpha$ is nontorsion.

\vskip .1 in 

Examples show that the conclusion of Theorem~\ref{MainThm} is false if
$\alpha$ is a torsion point, and that it can fail
if $\{\xi_n\}$ is merely a sequence of small points 
(that is, a sequence of points with $\hhat(\xi_n) \rightarrow 0$).  
In particular, Theorem~\ref{MainThm} 
cannot be strengthened to a theorem of Bogomolov type.  

\vskip .1 in

Theorem~\ref{MainThm} is the first known case of general conjectures 
by the second author (as refined by J. Silverman and S. Zhang) 
concerning abelian varieties and dynamical systems.
Assume as before that $k$ is a number field,
and let $S$ be a finite set of places of $k$ containing the archimedean places. 
Let $\cO_{k,S}$ be the ring of $S$-integers of $k$.   

\begin{conjecture} {\bf (Ih)}\label{Conj1}

Let $A/k$ be an abelian variety, and let  
$\cA_S/\Spec(\cO_{k,S})$ be a model of $A$.  
Let $D$ be an effective divisor on $A$, defined
over $\kbar$, at least one of whose irreducible components is not
the translate of an abelian subvariety by a torsion point, and 
let $\overline{D}$ be its Zariski closure in $\cA_S$.  
Then the set $A_{D,S}(\Zbar)_{\tors}$,  
consisting of all torsion points of $A(\kbar)$ whose 
closure in $\cA_S$ is disjoint from $\overline{D}$, 
is not Zariski dense in $A$.  
\end{conjecture} 

\begin{conjecture} {\bf (Ih)}\label{Conj2}

Let $R(x) \in k(x)$ be a rational function of degree
at least $2$, and consider the dynamical system associated to the 
rational map $R_* : \PP^1 \rightarrow \PP^1$.  
Let $\alpha \in \PP^1(\kbar)$ be non-preperiodic for $R_*$.
Then there are only finitely many pre-periodic points $\xi \in \PP^1(\kbar)$
which are $S$-integral with respect to $\alpha$, i.e. whose Zariski 
closures in $\PP^1/\Spec(\cO_{k,S})$ 
do not meet the Zariski closure of $\alpha$.
\end{conjecture} 

Theorem~\ref{MainThm}, 
in addition to being the one-dimensional case of Conjecture \ref{Conj1}, 
is equivalent to Conjecture \ref{Conj2} for Latt\`es maps.  
That is, if $E/k$ is an elliptic curve,  
let $R \in k(x)$ be the degree $4$ map 
on the $x$-coordinate corresponding to the doubling map on $E$, 
so that the following diagram commutes:
$$ 
  \begin{array}{ccc} 
       E        & \stackrel{[2]}{\rightarrow}  &  E  \\
   x \ \downarrow  \ \ \  &                              & \ \ \ \downarrow \ x \\ 
      \PP^1     & \stackrel{R_*}{\rightarrow}  &  \PP^1                
  \end{array} 
$$
Then $\beta \in E(\kbar)$ is a torsion point if and only $x(\beta)$ 
is preperiodic for $R_*$.

We will also prove Conjecture \ref{Conj2} for the map $R(x) = x^2$, 
where the preperiodic points are $0$, $\infty$, and the roots of unity.  
In that case, the assertion is that there are only finitely many roots of unity
which are $S$-integral with respect to a given non-root of unity 
$\alpha \in \kbar^{\times}$ (see Theorem~\ref{Thm1}).   
Conjecture \ref{Conj2} for Chebyshev maps
can be deduced by similar methods, though we do not do so here.  

\vskip .1 in
The motivation for Conjecture \ref{Conj1} is the following
analogy between diophantine theorems over $k$ and $\kbar$,
and over $\cO_k$ and $\Zbar$ (the ring of all algebraic integers).  
Let $A/k$ be an abelian variety, and let $X$ be a non-torsion subvariety of $A$ (that is, $X$ is not the translate of an abelian
subvariety by a torsion point).   
Recall that the Mordell-Lang Conjecture (proved by Faltings) says that 
$A(k) \cap X$ is not Zariski dense in $X$;  
while the Manin-Mumford Conjecture (first proved by Raynaud)
says that $A(\kbar)_{\tors} \cap X$ is not Zariski dense in $X$.  
Likewise, Lang's conjecture (also proved by Faltings) says that if $D$ is 
an effective ample divisor on $A$, 
then the set $A_D(\cO_k)$ of $\cO_k$-integral
points of $A$ not meeting $\supp(D)$ is finite. 
Note that $A$ is compact, whereas $A_D = A \backslash \supp(D)$ is noncompact. 

$$
\begin{tabular}{|l||c|c|} \hline
Type of variety; \  Type of rationality & $k$ & $\kbar$ \\ 
\hline 
Compact;  \qquad \ \ \ $k,\kbar$-rationality & Mordell-Lang &  Manin-Mumford \\ 
                             & Conjecture & Conjecture   \\ 
\hline  
Noncompact; \ \ \ \ $\cO_k,\Zbar$-rationality
       & Lang's     & Ih's   \\ 
       & Conjecture & Conjecture \ref{Conj1}  \\
\hline    
\end{tabular}
$$


\vskip .1 in

Conjecture~\ref{Conj2} is motivated by Conjecture~\ref{Conj1} and the familiar
analogy between torsion points of abelian varieties and preperiodic points of
rational maps.

\vskip .1 in
The paper is divided into two sections.  In the first, we prove  
Conjecture~\ref{Conj2} for the dynamical system  $R(x) = x^2$.  
In the second, we prove
Conjecture~\ref{Conj1} for elliptic curves.

Throughout the paper, we will use the following notation.  
For each place $v$ of $k$, let $k_v$ be the completion
of $k$ at $v$ and let $|x|_v$ be the normalized absolute value  
which coincides with the modulus of additive Haar measure on $k_v$.
If $v$ is archimedean and $k_v \cong \RR$, then $|x|_v = |x|$, 
while if $k_v \cong \CC$ then $|x|_v = |x|^2$.  If $v$ is nonarchimedean and lies
over the rational prime $p$, then $|p|_v = p^{-[k_v:\QQ_p]}$.  
For $0 \ne \alpha \in k$, the product formula reads 
\begin{equation*}
\prod_v |\alpha|_v \ = \ 1 \ .
\end{equation*} 
 
If $\kbar_v$ is an algebraic closure of $k_v$, there is
a unique extension of $|x|_v$ to $\kbar_v$, also denoted $|x|_v$. 
Given a finite extension $L/k$, for each place $w$ of $L$ we have the normalized 
absolute value $|x|_w$ on $L_w$.  
If we embed $L_w$ in $\kbar_v$, 
then $|x|_w = |x|_v^{[L_w:k_v]}$ for each $x \in L_w$. 
Write $\log(x)$ for the natural logarithm of $x$.
Given $\beta \in L$ and a place $v$ of $k$, 
as $\sigma$ ranges over all embeddings of $L$ into $\kbar_v$ fixing $k$ we have
\begin{equation} \label{FLH1} 
\sum_{\sigma : L/k \hookrightarrow \kbar_v} \log(|\sigma(\beta)|_v) 
\ = \ \sum_{w|v} \log(|\beta|_w) \ .
\end{equation}   

The absolute Weil height of $\alpha \in k$ (also called the naive height) is 
defined to be
\begin{equation*}
h(\alpha) \ = \ \frac{1}{[k:\QQ]} \sum_v \max(0, \log(|\alpha|_v)) \ . 
\end{equation*}
It is well known that for $\alpha \in \Qbar$, 
$h(\alpha)$ is independent of the field $k$ containing $\QQ(\alpha)$ 
used to compute it, so $h$ extends to a function on $\Qbar$.  
Furthermore $h(\alpha) \ge 0$, with $h(\alpha) = 0$ if and only if $\alpha = 0$ or $\alpha$ is a root of unity.

\vskip .1 in

\section{\bf Ih's conjecture for the dynamical system $R(x)=x^2$.} 

\subsection{\bf The finiteness theorem} 
\vskip .1 in

Let $S$ be a finite set of places of $k$ containing the archimedean places.
Given $\alpha, \beta \in \kbar$, view them as points in $\PP^1(\kbar)$ 
and let $\cl(\alpha), \cl(\beta)$ 
be their Zariski closures in $\PP^1/\Spec(\cO_k)$.  
By definition, $\beta$ is $S$-integral relative to $\alpha$ if 
$\cl(\beta)$ does not meet $\cl(\alpha)$ outside $S$.  
Thus, $\beta$ is $S$-integral relative to $\alpha$ if and only if
for each place $v$ of $k$ not in $S$, 
and each pair of embeddings $\sigma : k(\beta) \hookrightarrow \kbar_v$,
$\tau : k(\alpha) \hookrightarrow \kbar_v$, 
we have $\|\sigma(\beta),\tau(\alpha)\|_v = 1$ under the spherical metric on 
$\PP^1(\kbar_v)$. 
Equivalently, for all $\sigma$, $\tau$ 
\begin{equation*}
\left\{ \begin{array}{ll} 
 |\sigma(\beta)-\tau(\alpha)|_v \ge 1 & \text{if $|\tau(\alpha)|_v \le 1$\ ,} \\
 |\sigma(\beta)|_v \le 1 & \text{if $|\tau(\alpha)|_v > 1$\ .} 
        \end{array} \right.
\end{equation*}
\vskip .05 in
 
\begin{theorem} \label{Thm1}  
Let $k$ be a number field, and let $S$ be a finite set of places
of $k$ containing all the archimedean places.  
Fix $\alpha \in \kbar$ with $h(\alpha) > 0$;  
that is, $\alpha$ is not $0$ or a root of unity.  
Then there are only finitely many roots of
unity in $\kbar$ which are $S$-integral with respect to $\alpha$.  
\end{theorem} 

\vskip .05 in
Before giving the proof, we note some examples which
limit possible strengthenings of the theorem.

\vskip .1 in
A) \ The hypothesis $h(\alpha)>0$ is necessary:     

If $\alpha = 0$, take $k = \QQ$.  Then each root of unity $\zeta_n$ is integral
with respect to $\alpha$ at all finite places.  
If $\alpha = 1$, then each root of unity of composite order is integral with respect to $\alpha$ at all finite places. 

If $\alpha = \zeta_N$ is an $N^{th}$ root of unity with $N > 1$, 
take $k = \QQ(\zeta_N$).  
If $\zeta_m$ is a primitive $m^{th}$ root of unity with $(m,N) = 1$ and $m > 1$, 
then $\zeta_N^{-1}\zeta_m$ is a primitive $mN^{th}$ 
root of unity whose order divisible by at least two primes.  
This means $1-\zeta_N^{-1}\zeta_m$ is a unit, so $\zeta_N-\zeta_m$ is also a unit.  
Hence $\zeta_m$ is integral with respect to $\alpha$ 
at all finite places.

\vskip .1 in
B) \ When $h(\alpha) > 0$, one can ask if the theorem could be strengthened 
to a result of Bogomolov type:  is there a number $B = B(\alpha) > 0$ 
such that there are only finitely many 
points $\beta \in \kbar$ with $h(\beta) < B$ which are $S$-integral
with respect to $\alpha$?  That is, could  
finiteness for roots of unity be strengthened 
to finiteness for small points?  

\vskip .05 in
The following example shows this is not possible.  Take $k = \QQ$, 
$\alpha = 2$, and $S = \{\infty\}$.  
For each $n$, let $\beta_n$ be a root of the polynomial
\begin{equation*}
f_n(x)  \ =  \ x^{2^n}(x-2)-1  \ .
\end{equation*}
Here $f_n(x+1)$ is Eisenstein 
with respect to the prime $p=2$, so $f_n(x)$ is irreducible over $\QQ$. 
Note that each $\beta_n$ is a unit.  
By Rouch\'e's theorem, $\beta_n$ has one conjugate very
near $2$ and the rest of its conjugates very close to the unit circle;     
this can be used to show that $\lim_{n \rightarrow \infty} h(\beta_n) =0$.
Finally, $\beta_n - 2$ is also a unit, 
so $\beta_n$ is integral with respect to $2$
at all finite places.  

\vskip .1 in

\begin{proof} of Theorem \ref{Thm1}.

By replacing $k$ with $k(\alpha)$, 
and $S$ with the set of places $S_{k(\alpha)}$ lying over $S$, 
we are reduced to proving the theorem when $\alpha \in k$.  
Indeed, if $\zeta$ is a root of unity which is $S$-integral with respect 
to $\alpha$ over $k$, then each $k$-conjugate of $\zeta$ is 
$S_{k(\alpha)}$-integral with respect to $\alpha$ over $k(\alpha)$.  

Suppose $\alpha \in k$, and that there are infinitely many distinct 
roots of unity $\{\zeta_n\}$ which are $S$-integral with respect to $\alpha$.  
For each $n$, we will evaluate the sum
\begin{equation} \label{MainF}
A_n \ = \ \frac{1}{[k(\zeta_n):\QQ]} 
  \sum_{\text{$v$ of $k$}} \sum_{\sigma : k(\zeta_n)/k \hookrightarrow \kbar_v} 
                                 \log(|\sigma(\zeta_n)-\alpha|_v) 
\end{equation}
in two different ways.  On the one hand, an application of the product formula
will show that each $A_n = 0$. 
On the other hand, by applying the integrality hypothesis, Baker's theorem on 
linear forms in logarithms, 
and a strong form of equidistribution for roots of unity, 
we will show that $\lim_{n \rightarrow \infty} A_n = h(\alpha) > 0$.  
This contradiction will give the desired result.  

\vskip .1 in
The details are as follows.  
First, using (\ref{FLH1}), formula (\ref{MainF}) can be rewritten as
\begin{equation*}
A_n \ = \ \frac{1}{[k(\zeta_n):\QQ]} \sum_{\text{$w$ of $k(\zeta_n)$}} 
                                 \log(|\zeta_n-\alpha|_w) \ .
\end{equation*}
Since $\alpha$ is not a root of unity, the product formula gives $A_n = 0$.  

Next, take $v \notin S$.  
If $|\alpha|_v > 1$ then by the ultrametric inequality, 
for each $\sigma : k(\zeta_n)/k \hookrightarrow \kbar_v$ we have 
$|\sigma(\zeta_n)-\alpha|_v = |\alpha|_v$.  
On the other hand, if $|\alpha|_v \le 1$,  
the integrality hypothesis gives $|\sigma(\zeta_n)-\alpha|_v = 1$.  It follows 
that for each $v \notin S$  
\begin{equation}
\frac{1}{[k(\zeta_n):\QQ]} 
   \sum_{\sigma : k(\zeta_n)/k \hookrightarrow \kbar_v} 
                                 \log(|\sigma(\zeta_n)-\alpha|_v) 
       \ = \ \frac{1}{[k:\QQ]} \max(0,\log(|\alpha|_v)),
\end{equation} 
so that
\begin{eqnarray}
A_n & = &  
  \sum_{v \in S} \frac{1}{[k(\zeta_n):\QQ]} 
          \sum_{\sigma : k(\zeta_n)/k \hookrightarrow \kbar_v}  
                           \log(|\sigma(\zeta_n)-\alpha|_v)  \label{MainF1} \\
    &   & \qquad \qquad \qquad                               
          \ + \ \frac{1}{[k:\QQ]} \sum_{v \notin S} \max(0,\log(|\alpha|_v)) \ . \notag
\end{eqnarray}
Now let $n \rightarrow \infty$ in (\ref{MainF1}). Since $S$ is finite,   
we can interchange the limit and the sum over $v \in S$, obtaining
\begin{eqnarray} 
0  & = & \sum_{v \in S} \big( \lim_{n \rightarrow \infty} \frac{1}{[k(\zeta_n):\QQ]}
       \sum_{\sigma : k(\zeta_n)/k \hookrightarrow \kbar_v}  
                             \log(|\sigma(\zeta_n)-\alpha|_v) \big) \notag \\                                    
  &   & \qquad \qquad 
  \ + \ \frac{1}{[k:\QQ]} \sum_{v \notin S} \max(0,\log(|\alpha|_v) \ . \label{MainF2}
\end{eqnarray}
We will now show that for each $v \in S$,  
\begin{equation} \label{FGG1}
\lim_{n \rightarrow \infty} \frac{1}{[k(\zeta_n):\QQ]}
       \sum_{\sigma : k(\zeta_n)/k \hookrightarrow \kbar_v}  
                             \log(|\sigma(\zeta_n)-\alpha|_v) \ = \  
\frac{1}{[k:\QQ]} \max(0,\log(|\alpha|_v) \ . 
\end{equation} 
Inserting this in (\ref{MainF2}) gives $h(\alpha) = 0$, a contradiction.  

\vskip .1 in
For each nonarchimedean $v \in S$, 
(\ref{FGG1}) is trivial if $|\alpha|_v > 1$ or $|\alpha|_v < 1$.
In the first case $|\sigma(\zeta_n)-\alpha|_v = |\alpha|_v$ for all $n$ and all 
$\sigma$, and in the second case $|\sigma(\zeta_n)-\alpha|_v = 1$
for all $n$ and all $\sigma$.  Hence we can
assume that $|\alpha|_v = 1$.  

\begin{lemma} \label{Lem1}
Let $v$ be nonarchimedean, and suppose $|\alpha|_v = 1$.  Then

$(A)$  There is a bound $M(\alpha) > 0$ such that $|\zeta - \alpha|_v \ge M(\alpha)$ 
for all roots of unity $\zeta \in \kbar_v$.    

$(B)$  For each $0 < r < 1$, there are only finitely many roots of unity 
         $\zeta \in \kbar_v$ with $|\zeta-\alpha|_v < r$.  
\end{lemma}  

\begin{proof}  Since $\alpha$ is not a root of unity, (A) follows immediately
from (B).  For (B), note that if $\zeta$ and $\zeta^{\prime}$ are roots of unity
with $|\zeta-\alpha|_v < r$ and $|\zeta^{\prime} - \alpha|_v < r$,  
then $|\zeta-\zeta^{\prime}|_v < r$ and so 
$\zeta^{\prime \prime} = \zeta^{-1} \zeta^{\prime}$ is a root of unity with  
$|1 - \zeta^{\prime \prime}|_v < r$.  There are only finitely 
many such $\zeta^{\prime \prime}$.  Indeed, if $p$ is the rational prime
under $v$, the only roots of unity $\xi \in \kbar_v$ with 
$|1-\xi|_v < 1$ are ones with order  $p^n$  for some  $n$.  If $\xi$ is a 
primitive  $p^n$-th root of unity, 
then $|1-\xi|_v = p^{-[k_v:\QQ_p]/p^{n-1}(p-1)}$
so $1 > r > |1-\xi|_v$ for only finitely many $n$.   
\end{proof}         
 
\vskip .1 in
Assuming $v$ is nonarchimedean and $|\alpha|_v = 1$, 
let $M(\alpha)$ be as in the Lemma.
Fix $0 < r < 1$, and let $N(r)$ be the number of roots of unity in $\kbar_v$
with $|\zeta-\alpha|_v < r$.  For each $\zeta_n$ 
and each $\sigma : k(\zeta_n)/k \rightarrow \kbar_v$, we have
$|\sigma(\zeta_n) - \alpha|_v \le 1$, so 
\begin{eqnarray*}
0 & \ge & \lim_{i \rightarrow \infty} \frac{1}{[k(\zeta_n:\QQ]}
       \sum_{\sigma : k(\zeta_n)/k \hookrightarrow \kbar_v}  
                             \log(|\sigma(\zeta_n)-\alpha|_v) \\
 & \ge & \lim_{n \rightarrow \infty} \frac{1}{[k(\zeta_n):\QQ]}
         (([k(\zeta_n):k]-N(r)) \cdot \log(r)+N(r)\cdot \log(M(\alpha))) \\
 & = & \frac{1}{[k:\QQ]} \log(r) \ .
\end{eqnarray*}   
Since $r < 1$ is arbitrary, the limit in (\ref{FGG1}) is $0$, verifying 
(\ref{FGG1}) in this case.

\vskip .1 in
Now suppose $v$ is archimedean.  To simplify notation, 
view $k$ as a subfield of $\CC$ and identify $\kbar_v$ with $\CC$.  
(Thus, the way $k$ is embedded depends on the choice of $v$).   

By Jensen's formula (see \cite{Conway}, p.280)
applied to $f(z) = z-\alpha$,  
\begin{equation} \label{FGG2} 
\frac{1}{2\pi} \int_0^{2\pi} \log(|e^{i\theta} - \alpha|) \, d\theta
\ = \ \max(0,\log(|\alpha|)) \ .  
\end{equation}
Here $|x|$ can be replaced by $|x|_v$.  

  The $\Gal(\kbar/k)$-conjugates of roots of unity 
equidistribute in the unit circle.  
We will give a direct proof of this below, but we note that it also follows
from generalizations of Bilu's theorem,   
for example the equidistribution theorem for polynomial dynamical systems 
given in Baker-Hsia (\cite{B-H}).  
The Baker-Hsia theorem implies that if $\mu_n$ is the discrete measure
\begin{equation*}
\mu_n  \ = \ \frac{1}{[k(\zeta_n):k]} 
\sum_{\sigma: k(\zeta_n)/k \hookrightarrow \CC} \delta_{\sigma(\zeta_n)}(x) \ , 
\end{equation*}
where $\delta_P(x)$ is the Dirac measure with mass $1$ at $P$, then the $\mu_n$
converge weakly to the Haar measure $\mu = (1/2\pi) d\theta$ on the unit circle.

If $|\alpha|_v > 1$ or $|\alpha|_v < 1$ then  $\log(|z-\alpha|_v)$ 
is continuous on the unit circle.  In these cases, (\ref{FGG1}) follows
from (\ref{FGG2}) and weak convergence.  If $|\alpha|_v = 1$ then 
$\log(|z-\alpha|_v)$ is not continuous on $|z|=1$ and weak convergence
is not enough to give  
$\int_{|z|=1} \log(|z-\alpha|) \, d\mu_n(z) \rightarrow 0$.  
There could be a problem if some conjugate 
were extremely close to $\alpha$, 
or if many conjugates clustered near $\alpha$.  

\vskip .1 in
The first problem is solved by A. Baker's theorem on lower bounds for
linear forms in logarithms (see Baker \cite{Baker}, Theorem 3.1, p.22). 
We are assuming that $|\alpha|_v = 1$, and $\alpha$ is not a 
root of unity.  Fix a branch of $\log$
with $\log(z) = \log(|z|)+i\theta$, $-\pi < \theta \le \pi$, 
and write $\log(\alpha)=i\theta_0$.  For another branch, $\log(1)=2\pi i$.
The following is a special case of Baker's theorem.  
(In his statement of the theorem, Baker uses an exponential height 
having bounded ratio with $H(\beta) = e^{h(\beta)}$.) 

\begin{proposition} \label{Prop1} \text{\rm (A. Baker)}
There is a constant $C = C(\alpha) > 0$  
such that for each $\beta = a/N \in \QQ$  
\begin{equation*}
|\, i \theta_0 - \beta \cdot 2 \pi i|  \ \ge \ 
 e^{-C \cdot \max(1, h(\beta))} \ ,  
\end{equation*}
where $h(\beta) = \log(\max(|a|,|N|))$ is the absolute height of $\beta$.
\end{proposition} 

\vskip .1 in
The second problem is settled by 
a strong form of equidistribution for roots of unity, 
proved in \S \ref{Sec1.2} below.  
It says that for any $0 < \gamma < 1$, 
the conjugates of the $\zeta_n$ are asymptotically 
equidistributed in arcs of length $[k(\zeta_n):k]^{-\gamma}$. 
Note that weak convergence is equivalent to equidistribution  
in arcs of fixed length.   

\begin{proposition}  \label{Prop2} {\rm (Strong Equidistribution)}  
Let $k \subset \CC$ be a number field.  
Then the $\Gal(\kbar/k)$-conjugates of the roots of unity in $\kbar$
$($viewed as embedded in $\CC)$
are strongly equidistributed in the unit circle, in the following sense.  

Given an arc $I$ in the unit circle, write 
$\mu(I)=\frac{1}{2\pi}\text{\rm length}(I)$ for its normalized Haar measure.
If $\zeta \in \kbar$ is a root of unity, put 
\begin{equation*}
N(\zeta,I) \ = \ \#\{\sigma(\zeta) \in I : \sigma \in \Gal(\kbar/k) \}  \ .
\end{equation*}
Fix $0 < \gamma < 1$.  Then for all roots of unity $\zeta$ and all $I$, 
\begin{equation} \label{FMJ2B}
\frac{N(\zeta,I)}{[k(\zeta):k]} \ = \  
            \mu(I) + O_{\gamma}([k(\zeta):k]^{-\gamma}) \ . 
\end{equation}
\end{proposition}  

\vskip .1 in
We remark that a strong equidistribution theorem for points of small height with respect to an 
arbitrary dynamical systems on $\PP^1$ 
has recently been proved by C. Favre and J. Rivera-Letelier 
(\cite{F-RL}, Th\'eor\`eme 6).   

\vskip .1 in
Assuming Proposition~\ref{Prop2}, we will now complete the proof of Theorem \ref{Thm1}
by showing that (\ref{FGG1}) holds for archimedean $v$ when $|\alpha|_v=1$. 
 
Let $\mu = (1/2\pi) d \theta$ be the normalized Haar measure on the unit circle,
and for each $n$, put 
\begin{equation*}
\mu_n \ = \ \frac{1}{[k(\zeta_n):k]} 
     \sum_{\sigma : k(\zeta_n)/k \rightarrow \CC} \delta_{\sigma(\zeta_n)}(x) \ .
\end{equation*}
Then $\mu_n$ is supported on the unit circle and the $\mu_n$ converge
weakly to $\mu$. We must show that 
\begin{equation*}
\int_{|z|=1} \log(|z-\alpha|) \, d\mu_n(z) \ = \
   \frac{1}{[k(\zeta_n):k]} \sum_{\sigma} \log(|\sigma(\zeta_n)-\alpha|) 
      \ \rightarrow \ 0 \ .  
\end{equation*}

The idea is to break the sum into three parts:  the terms nearest $\alpha$, 
which can be treated by Baker's theorem;  the other terms in a small neighborhood
of $\alpha$, which can be dealt with by strong equidistribution;  and the
rest, which can be handled by weak convergence. 
 
Fix $0 < \epsilon < 1$.  We will show that for all sufficiently large $n$,
\begin{equation} \label{FGL4}
|\int_{|z|=1} \log(|z-\alpha|) \, d\mu_n(z)|  \ < \ 6 \epsilon \ . 
\end{equation}

Note that $\int_0^{\varepsilon} \log(t/\varepsilon) \, dt =  -\varepsilon$. 
For the remainder of the proof, we restrict to $|z| = 1$;  
write $\alpha = e^{i \theta_0}$ where $-\pi < \theta_0 \le \pi$, 
and write $z = e^{i \theta}$ where $\theta_0 - \pi < \theta \le \theta_0 + \pi$. 
Define   
\begin{equation*}
\larg_{\alpha,\varepsilon}(z) \ = \min(0,\log(|\theta-\theta_0|/\varepsilon)) \ .
\end{equation*} 
Then there is a continuous function $g_{\alpha,\varepsilon}(z)$ on $|z| = 1$ for 
which 
$\log(|z-\alpha|) = \larg_{\alpha,\varepsilon}(z) + g_{\alpha,\varepsilon}(z) $.   
Recalling that $\int_{|z|=1} \log(|z-\alpha|) \, d\mu(z) = 0$, we have  
\begin{equation*}
\int_{|z| = 1} g_{\alpha,\varepsilon}(z) \, d\mu(z) \ = \ 
- \int_{|z| = 1} \larg_{\alpha,\varepsilon}(z) \, d\mu(z) \ = \ 
-2 \int_0^{\varepsilon} \log(\theta/\varepsilon) \, \frac{d \theta}{2 \pi} 
\ = \ \frac{\varepsilon}{\pi} \ .
\end{equation*} 
By weak convergence, it follows that for all sufficiently large $n$, 
\begin{equation} \label{FMP1} 
|\int_{|z| = 1} g_{\alpha,\varepsilon}(z) \, d\mu_n(z)| \ < \ \varepsilon \ .
\end{equation}

To obtain (\ref{FGL4}), it will suffice to show that for all sufficiently large $n$,
\begin{equation*} 
|\int_{|z| = 1} \larg_{\alpha,\varepsilon}(z) \, d\mu_n(z)| \ < \ 5 \varepsilon \ .
\end{equation*}
For each interval $[c,d]$ let $I_{\alpha}([c,d])$ be the arc 
$\{\alpha e^{2 \pi i t} : t \in [c,d] \}$.   
Noting that $\larg_{\alpha,\varepsilon}(z)$ is supported on
$I_{\alpha}([-\varepsilon,\varepsilon])$, 
put $D = D_n = \lceil [k(\zeta_n):k]^{1/2} \rceil$ 
and divide $I_{\alpha}([-\varepsilon,\varepsilon])$ into $2D$ equal subarcs.  
Taking $\gamma = 2/3$ in Proposition \ref{Prop2}, 
it follows that if $n$ is sufficiently large, each such subarc contains at most
$2\varepsilon[k(\zeta_n):k]^{1/2}$ conjugates of $\zeta_n$.  

First consider the union of the two central subarcs, 
$I_{\alpha}([-\varepsilon/D,\varepsilon/D])$.  Let
$N$ be the order of $\zeta_n$.  Let $\sigma_0(\zeta_n) = e^{2\pi i a/N}$ 
be the conjugate of $\zeta_n$ closest to $\alpha = e^{i \theta_0}$.  
We can assume that $|a/N| \le 1$, which implies that   
that $h(a/N) = \max(\log(|a|),\log(N)) = \log(N)$.  By Baker's theorem, 
\begin{equation*}
|2\pi (a/N) - \theta_0| \ > \  e^{-C \max(1, \log(N))} \ .
\end{equation*}
Hence if $n$ is sufficiently large,  
\begin{equation*}
\larg_{\alpha,\varepsilon}(\sigma_0(\zeta_n)) 
\ > \ - C \log(N) - \log(\varepsilon) 
\ \ge  \ -C \log(N) \ .
\end{equation*}                               
Since there are at most $4 \varepsilon [k(\zeta_n):k]^{1/2}$ conjugates of $\zeta_n$
in $I_{\alpha}([-\varepsilon/D,\varepsilon/D])$,  
\begin{equation*} 
0 \ \ge \ \int_{I_{\alpha}([-\varepsilon/D,\varepsilon/D])} 
\larg_{\alpha,\varepsilon}(|z-\alpha|) \, d\mu_n(z)
\ > \ -4 \frac{ C \log(N) }{[k(\zeta_n):k]^{1/2}} \varepsilon  \ .
\end{equation*}
Note that $[k(\zeta_n):k] \ge [\QQ(\zeta_n):\QQ]/[k:\QQ] = \varphi(N)/[k:\QQ]$.
For all large $N$, $\varphi(N) \ge N^{1/2}$, so there
is a constant $B$ such that $[k(\zeta_n):k]^{1/2} \ge B N^{1/4}$.  
Thus for all sufficiently large $n$, 
\begin{equation} \label{FGL5}
| \int_{I_{\alpha}([-\delta/D,\delta/D])} \log(|z-\alpha|) \, d\mu_n(z)|
\ < \ \epsilon \ .
\end{equation}

Finally, consider the remaining subarcs.
For $\ell = 1, \ldots, D-1$, if 
\begin{equation*}
z \in I_{\alpha}([\ell \varepsilon/D,(\ell+1) \varepsilon/D]) \quad \text{\rm or} \quad  
z \in I_{\alpha}([-\ell \varepsilon/D,-(\ell+1) \varepsilon/D]) 
\end{equation*}
then $0 \ge \larg_{\alpha,\varepsilon}(z) \ge \log(\ell/D)$.
As before, by Proposition \ref{Prop2}, 
for sufficiently large $n$, each subarc contains at most 
$2 [k(\zeta_n):k] (\varepsilon/D)$ conjugates of $\zeta_n$.  
It follows that 
\begin{eqnarray}
0 & \ge & \int_{I_{\alpha}([-\delta,\delta]) \backslash 
      I_{\alpha}([-\delta/D,\delta/D])} \larg_{\alpha,\varepsilon}(z) \, d\mu_n(z)
                      \notag \\
  & \ge & 2 \cdot \sum_{\ell=1}^{D-1} 
      \log((\frac{\ell \varepsilon}{D})/\varepsilon) \cdot \frac{2 \varepsilon}{D} 
                      \notag \\
  & > & 4 \int_0^{\varepsilon} \log(t/\varepsilon) \, dt  \ = \ - 4 \epsilon \ . 
                      \label{FGL6}
\end{eqnarray}                       
Combining (\ref{FMP1}), (\ref{FGL5}), and (\ref{FGL6}) gives (\ref{FGL4}), 
which completes the proof of Theorem~\ref{Thm1}.
\end{proof}

\vskip .1 in
In the course of writing this paper, 
the authors learned of several results related to Theorem \ref{Thm1},
some of which imply it in special cases.  

A.~Bang's theorem \cite{Bang} (1886) 
says that if $\alpha \ne \pm 1$ is a nonzero rational number, 
then for all sufficiently large integers $n$ there is a prime $p$ 
such that the order of $\alpha$ modulo $p$ is exactly $n$.  
This can be rephrased as saying that for all 
sufficiently large $n$, there exists a primitive $n$-th root of unity $\zeta_n$
and a nonzero prime ideal $\fp$ of $\ZZ[\zeta_n]$ such that $\alpha \equiv \zeta_n \pmod{\fp}$.  
Since all primitive $n$-th roots are conjugate over $\QQ$, 
this implies Theorem \ref{Thm1} in the case
$\alpha \in \QQ$.  A.~Schinzel \cite{Schinzel} gave an effective  
generalization of Bang's theorem to arbitrary number fields;  Schinzel's theorem 
implies Theorem \ref{Thm1} for number fields $k$ which are linearly disjoint
from the maximal cyclotomic field $\QQ^{ab}$, and $\alpha \in k$.  

J. Silverman \cite{Sil3} has shown that if $\alpha \in \Qbar$ 
is an algebraic unit which is not a root of unity, there are only finitely
many $m$ for which $\Phi_m(\alpha)$ is a unit, where $\Phi_m(x)$ is the $m$-th
cyclotomic polynomial.  In fact, if $d = [\QQ(\alpha):\QQ]$ he shows there
is an absolute, effectively computable constant $C$ such that the number
of such $m$'s is at most 
\begin{equation*}
C \cdot d^{1 + 0.7/\log(\log(d))} \ .
\end{equation*}
In the case when $\alpha$ is a unit, this yields Theorem~\ref{Thm1} in the same situations as Schinzel's theorem.

G. Everest and T. Ward (\cite{E-W}, Lemma 1.10)
show that if $F(x) \in \ZZ[x]$ is monic and 
irreducible, with roots ${\alpha_1, \ldots, \alpha_d}$, 
and if $F(x)$ is not a constant multiple of 
$x$ or a cyclotomic polynomial
$\Phi_m(x)$, then the quantity
$\Delta_n(F) = \prod_{i=1}^d (\alpha_i^n-1)$ satisfies
\begin{equation} \label{BB2}
\lim_{n \rightarrow \infty} \frac{1}{n} \log(\Delta_n(F)) \ = \ m(F) \ > \ 0,
\end{equation}
where $m(F) = \deg(F) \cdot h(\alpha_i)$ is the logarithm of the 
Mahler measure of $F(x)$.   
When $k = \QQ$, and $\alpha = \alpha_1$ is an algebraic integer, 
the product formula tells us that 
$\prod_{\text{$v$ of $\QQ$}} |\Delta_n(F))|_v = 1$, so for all large $n$ there 
must be some nonarchimedean $v$ and some $\alpha_i$ such that 
$|\alpha_i^n-1|_v = 1$, and this in turn means there is some 
$n$-th root of unity $\zeta$ with $|\alpha_i-\zeta|_v < 1$.  
However, this is not strong enough to give Theorem \ref{Thm1} because
(a) $\zeta$ might not be primitive, 
and (b) the primitive $n$-th roots of unity might
not all be conjugate to one another over $\QQ(\alpha)$.
  
\vskip .1 in
\subsection{\bf   Strong equidistribution for roots of unity.} \label{Sec1.2}

We will now prove Proposition~\ref{Prop2}, 
the strong equidistribution theorem for roots of unity. 
At least when $k = \QQ$, the result is well known to analytic number theorists, 
but we do not know a reference in the literature.  

The proof rests on the following lemma, for which we thank Carl Pomerance. 
Let $\varphi(N)$ denote Euler's function 
and let $d(N) = \sum_{m|N} 1$ be the divisor function.  We write 
$\lambda(m)$ for the number of primes dividing $m$, and use $\theta(x)$
to denote a quantity satisfying $-x \le \theta(x) \le x$. 
 
\begin{lemma} \label{Lem2} {\rm (Pomerance)}
Fix an integer $Q > 1$ and an integer $b$ coprime to $Q$.     
Then for each integer $N \ge 1$ divisible by $Q$ 
and each interval $(c,d] \subset \RR$,
\begin{equation*} 
\#\{a \in (c,d] \cap \ZZ :  (a,N) = 1, a \equiv b \ \text{\rm{(mod $Q$)}} \}
 \ = \ \frac{\varphi(N)}{N \varphi(Q)} (d-c) + \theta(d(N))  \ .
\end{equation*} 
\end{lemma}

\begin{remark}
The main content of the lemma is that the error depends only on $N$, 
and not on $Q$ or $(c,d]$.  
\end{remark}

\begin{proof}  
 Let  $p_1, \ldots, p_r$ be the primes dividing  $N$ but not $Q$. 
(If there are no such primes,  
take $p_1 \cdots p_r = 1$ in the argument below).   
Take $b_0 \in \ZZ$ with  $b_0 \equiv b \text{\ (mod $Q$)}$, 
$b_0 \equiv 0 \text{\ (mod $p_1 \ldots p_r$)}$.  Then
\begin{eqnarray*}
& & \{a \in (c,d] \cap \ZZ : a \equiv b \text{\ (mod $Q$)}, (a,N) = 1\} \\
& & \qquad \qquad \qquad \ = \ \{a \in (c,d] \cap \ZZ :\ Q|a-b_0, \  
p_1, \ldots, p_r \ndiv a-b_0 \, \}
\end{eqnarray*}
For each positive integer $m$ dividing $p_1 \cdots p_r$ put 
$r_{m,b,Q}(c,d) = \#\{a \in (c,d] \cap \ZZ :  Qm|a-b_0 \}$.  Then 
\begin{equation*}
r_{m,b,Q}(c,d) \ = \ \lfloor (d-b_0)/Qm \rfloor - \lfloor (c-b_0)/Qm \rfloor 
     \ = \ \frac{1}{Qm}(d-c) + \theta(1) \ .
\end{equation*} 
Carrying out inclusion/exclusion relative to the primes $p_1, \ldots, p_r$ 
we have  
\begin{eqnarray*}
& & \#\{a \in (c,d] \cap \ZZ : a \equiv b \text{\ (mod $Q$)}, (a,N) = 1\} \\
& & \qquad \qquad \qquad \qquad \qquad
       \ = \ \sum_{m | p_1 \cdots p_r} (-1)^{\lambda(m)} \, r_{m,b,Q}(c,d) \\
& & \qquad \qquad \qquad \qquad \qquad
          \ = \ \frac{1}{Q} \prod_{i=1}^r (1-\frac{1}{p_i}) \cdot (d-c)  
                                     + \theta(d(p_1 \cdots p_r)) \\
& & \qquad \qquad \qquad \qquad \qquad
          \ = \  \frac{\varphi(N)}{N \varphi(Q)} (d-c) 
                                     + \theta(d(N)) \ .
\end{eqnarray*}                                                                          
\end{proof} 
  
\begin{proof} of Proposition \ref{Prop2}.  

Let $\zeta_N$ denote a primitive $N^{th}$ root of unity. 
There are only finitely many subfields of $k$, so there are 
only finitely subfields of the form 
$k_N = k \cap \QQ(\zeta_N)$ for some $N$.  For each $N$ there 
is a minimal $Q$ for which $k_N = k_Q$, and then 
$\QQ(\zeta_Q) \subset \QQ(\zeta_N)$ so $Q|N$. We will call $Q = Q_N$
the cyclotomic conductor of $\zeta_N$ relative to $k$, and write
$T_N = [\QQ(\zeta_{Q_N}):k_N]$.  

As $\QQ(\zeta_N)$ is Galois over $\QQ$, it is linearly disjoint from 
$k$ over $k_N$, and $\Gal(k(\zeta_N)/k) \cong \Gal(\QQ(\zeta_N)/k_N)$.  
Since $k_N \subset \QQ(\zeta_{Q_N}) \subset \QQ(\zeta_N)$, 
the conjugates of $\zeta_N$ over $k$ are a union of $T_N$ sets of the form 
\begin{equation*}
\{e^{2\pi i a/N} : a \equiv b_i \text{\rm (mod $Q_N$)}, (a,N) = 1 \} \ ,   
\end{equation*}
for certain numbers $b_i$ coprime to $Q_N$.

Let $I$ be an arc of the unit circle corresponding to an 
angular interval $(\theta_1,\theta_2]$.  
Put $(c,d] = \frac{N}{2\pi}(\theta_1,\theta_2]$.  
Then $e^{2 \pi i a/N} \in I$ if and only if $a \in (c,d]$. 
By Lemma \ref{Lem2}, 
\begin{equation} \label{FMJ1}
N(\zeta_N,I) \ = \ 
T_N \cdot \frac{\varphi(N)}{N \varphi(Q_N)} 
        \cdot \frac{N}{2 \pi}(\theta_2-\theta_1) 
\ + \  \theta(T_N \cdot d(N)) \ .
\end{equation}

Recall that for any $\delta>0$, if $N$ is sufficiently large then   
$d(N) \le N^{\delta}$ and $\varphi(N) \ge N^{1-\delta}$  
(see Hardy and Wright \cite{H-W}, Theorem 315, p.260, and Theorem 327, p.267).
Take $\delta$ such that  $0 < 2\delta < 1-\gamma$.   
Noting that $[k(\zeta_N):k]=T_N \varphi(N)/\varphi(Q_N)$, 
and that $\varphi(Q_N)$ is bounded independent of $N$, 
(\ref{FMJ1}) gives 
\begin{equation} \label{FMJ2}
\frac{N(\zeta_N,I)}{[k(\zeta_N):k]} \ = \ \mu(I) + O_{\gamma}(N^{-\gamma}) \ .
\end{equation}
Since $[k(\zeta_N):k]\le N$, 
the error bound in (\ref{FMJ2}) holds with $N$ replaced by $[k(\zeta_N):k]$. 
Since $[k(\zeta_N):k]/N^{\gamma} \rightarrow \infty$ as $N \rightarrow \infty$, 
adjoining or removing endpoints of $I$ will not affect the form of the 
estimate, so (\ref{FMJ2B}) applies to all intervals.  
\end{proof}

\section{\bf  Ih's conjecture for elliptic curves.}
\subsection{\bf The finiteness theorem} 

\vskip .05 in
Let $k$ be a number field, and let $E/k$ be an elliptic curve.  
We can assume $E$ is defined by a Weierstrass equation
\begin{equation} \label{FE1} 
y^2 + a_1 xy + a_3 y \ = \ x^3 + a_2 x^2 + a_4 x + a_6
\end{equation}
with coefficients in $\cO_k$.  
More precisely, $E$ is the hypersurface in $\PP^2/\Spec(k)$ defined by the 
homogenization of (\ref{FE1}). 
Let $\Delta$ be its discriminant.   

Given a nonarchimedean place $v$ of $k$ and points 
$\alpha, \beta \in E(\kbar)$, 
we will say that $\beta$ is integral with respect to $\alpha$ at $v$ if
the Zariski closures $\cl(\beta)$ and $\cl(\alpha)$ do not meet in the model 
$\cE_v/\Spec(\cO_v)$ defined by the homogenization of (\ref{FE1}).  
Equivalently, if $\|z,w\|_v$ is the restriction 
of the spherical metric on $\PP^2(\kbar_v)$ to $E(\kbar_v)$ 
(see \cite{R1}, \S1.1), 
then for each pair of embeddings $\sigma, \tau : \kbar/k \hookrightarrow \kbar_v$,   
\begin{equation*}
\|\sigma(\beta),\tau(\alpha)\|_v \ = \ 1 \ .
\end{equation*}
If $S$ is a set of places of $k$ containing all the archimedean places, 
we say $\beta$ is $S$-integral with respect to $\alpha$ if $\beta$ 
is integral with respect to $\alpha$ at each $v \notin S$.  

Write $\hhat(\alpha)$ for the canonical height on $E(\kbar)$, defined by
\begin{equation*}
\hhat(\alpha) \ = \ \frac{1}{2} 
\lim_{n \rightarrow \infty} \frac{1}{4^n} h_{\PP^1}(x([2^n]\alpha)) 
\ = \ 
\frac{1}{3} \lim_{n \rightarrow \infty} \frac{1}{4^n} h_{\PP^2}([2^n]\alpha) \ ,
\end{equation*}
where $h_{\PP^1}$ (resp. $h_{\PP^2}$) 
is the naive height on $\PP^1(\kbar)$ (resp. $\PP^2(\kbar)$), 
and $[m]$ is multiplication by $m$ on $E(\kbar)$.  
(For a discussion of $\hhat(\alpha)$ and its properties, 
see \cite{Sil1}, pp.227-231 and 365-366; 
or \cite{Sil2}, \S VI.)  
Recall that $\hhat(\alpha) \ge 0$, 
that $\hhat([m]\alpha) = m^2\hhat(\alpha)$ for all $m$, and that 
$\hhat(\alpha) = 0$ if and only if $\alpha \in E(\kbar)_{\tors}$.  
From these facts it follows (as is well known) 
that if $\xi \in E(\kbar)_{\tors}$, then 
\begin{equation} \label{FCOM1} 
\hhat(\alpha) \ = \ \hhat(\alpha-\xi) \ .  
\end{equation}

There is also a decomposition of $\hhat(\alpha)$ as a sum of local terms.
For each place $v$ of $k$, let $\lambda_v(P)$ be the local
N\'eron-Tate height function on $E(\kbar_v)$. 
For compatibility with our absolute values  
we normalize $\lambda_v(P)$ so that 
$\lambda_v(P) = [k_v:\QQ_p] \cdot \lambda_{v,\text{Sil}}(P)$, 
where $\lambda_{v,\text{Sil}}(P)$ 
is the local N\'eron-Tate height defined in Silverman (\cite{Sil1}, p.365).
For each $0 \ne \alpha \in E(k)$
\begin{equation} \label{FCOM3}
\hhat(\alpha) \ = \ 
\frac{1}{[k:\QQ]} \sum_{\text{$v$ of $k$}} \lambda_v(\alpha)
\end{equation}
(see \cite{Sil1}, Theorem 18.2, p.365).  
Note that only finitely many terms in the sum are nonzero.  

If $L/k$ is a finite extension, for each place $w$ of $L$ there is a   
normalized local N\'eron-Tate height $\lambda_w(P)$ on $E(\Lbar_w)$.  
If we fix an isomorphism $\Lbar_w \cong \kbar_v$, then 
for all $P \in E(\kbar_v)$,   
\begin{equation} \label{FCOM4} 
\lambda_w(P) \ = \ [L_w:k_v] \lambda_v(P) \ .
\end{equation}  
It follows that if $\beta \in E(L)$, 
then for each place $v$ of $k$, as $\sigma$ runs over
all embeddings of $L$ into $\kbar_v$ fixing $k$, 
\begin{equation} \label{FAA1} 
\sum_{\sigma:L/k \hookrightarrow \kbar_v} \lambda_v(\sigma(\beta)) 
     \ = \ \sum_{w|v} \lambda_w(\beta) \ .  
\end{equation}

We will use the following explicit formulas.

\begin{proposition} \label{LocFormulas}

Let $k$ be a number field, and let $E/k$ be an elliptic curve. 
Let $v$ be a place of $k$.   

$A)$ If $v$ is archimedean, fix an 
an isomorphism $E(\kbar_v) \cong \CC/\Lambda$ for an appropriate 
lattice $\Lambda \subset \CC$.  Let $\sigma(z,\Lambda)$
be the Weierstrass $\sigma$-function, 
let $\Delta(\Lambda) = g_2(\Lambda)^3-27 g_3(\Lambda)^2$ be the discriminant 
of $\Lambda$, and let $\eta : \CC \rightarrow \RR$ be the $\RR$-linearized 
period map associated to the Weierstrass $\zeta$-function $\zeta(z,\Lambda)$.  
If $P \in E(\kbar_v)$ corresponds to $z \in \CC$, then
\begin{equation*}
\lambda_v(P) \ = \ 
-\log(|\Delta(\Lambda)^{1/12} e^{-z\eta(z)/2} \sigma(z,\Lambda)|_v) \ .
\end{equation*}

Furthermore, if $\mu_v(z)$ is the additive Haar measure on $E(\kbar_v)$ 
which gives $E(\kbar_v) \cong \CC/\Lambda$ total mass $1$, then 
\begin{equation*}
\int_{E(\kbar_v)} \lambda_v(z) \, d\mu_v(z) \ = \ 0 \ .
\end{equation*}

$B)$ If $v$ is nonarchimedean and $E$ has split multiplicative reduction
at $v$ $($so $E$ is $k_v$-isomorphic to a Tate curve$)$, 
fix a Tate isomorphism
$E(\kbar_v) \cong \kbar_v^{\times}/q^{\ZZ}$ where $q \in \kbar_v^{\times}$ 
satisfies  $|q|_v = |1/j(E)|_v < 1$.  Let $B_2(x) = x^2 - x + \frac{1}{6}$ 
be the second Bernoulli polynomial, and put 
$\tlambda_v(x) 
= \frac{1}{2} B_2(\frac{x}{\ord_v(q)}) \cdot (-\log(|q|_v))$.    
If $P \in E(\kbar_v)$ corresponds to $z \in \kbar_v^{\times}$,  
with $z$ chosen so that $|q|_v < |z|_v \le 1$, then 
\begin{equation*}  
\lambda_v(P) \ = \  -\log(|1-z|_v) + \tlambda_v(\ord_v(z)) \ .
\end{equation*}

Furthermore, if $\mu_v$ is the Haar measure $dx/\ord_v(q)$ 
giving the loop $\RR/(\ZZ \cdot \ord_v(q))$ total mass $1$, then
\begin{equation*}
\int_{0}^{\ord_v(q)} \tlambda_v(x) \, d\mu_v(x) \ = \ 0 \ .
\end{equation*}

$C)$ If $v$ is nonarchimedean and $E$ has good reduction at $v$, 
let $\|z,w\|_v$ be the spherical metric on $E(\kbar_v)$ induced by 
a projective embedding $E \hookrightarrow \PP^2$ corresponding to 
a minimal Weierstrass model for $E$ at $v$.   
Then for each $P \in E_v(\kbar_v)$
\begin{equation*} 
\lambda_v(P) \ = \ -\log(\|P,0\|_v) \ .  
\end{equation*}
\end{proposition}  

\begin{proof}  This is a summary of results in (\cite{Sil2}, \S VI);  
see in particular Theorem 1.1, p. 455;  Theorem 3.2, p.466;
Theorem 3.3, p.468;  and Theorem 4.1, p.470.  
\end{proof} 

\vskip .1 in
We now come to Ih's conjecture for elliptic curves.  The following is a restatement
of Theorem \ref{MainThm} in the Introduction.  

\begin{theorem} \label{Thm2}   Let $E/k$ be an elliptic curve, and
let $S$ be a finite set of places of $k$, containing all the archimedean
places.  Let $\alpha \in E(\kbar)$ be a nontorsion point, i.e., a point with 
$\hhat(\alpha) > 0$.  Then there are only finitely 
many $\xi \in E(\kbar)_{\tors}$ which are $S$-integral with
respect to $\alpha$.  
\end{theorem}

Again there are limitations to possible strengthenings of the theorem:  

\vskip .1 in
A) \ As noted by Silverman, it is necessary that  $\alpha$ be nontorsion.     
If $\alpha = 0$ and $S$ is the set of archimedean places, 
then by Cassels' generalization of the Lutz-Nagell theorem
(Proposition \ref{CasselsThm} below), 
each torsion point whose order is divisible by at least two primes is 
$S$-integral with respect to $\alpha$.  

Similarly, if $\alpha$ is a torsion point of order $N > 1$, 
let $S$ contain all places of bad reduction for $E$.    
Then for each $q$ coprime to $N$, all $q$-torsion points are $S$-integral 
with respect to $\alpha$.  

\vskip .1 in
B) \ When $\hhat(\alpha) > 0$, Zhang has pointed out 
that Theorem \ref{Thm2} cannot be strengthened 
to a result of Bogomolov type.
A result of E. Ullmo (\cite{Ull0}, Theorem 2.4) 
shows that for each $\varepsilon > 0$,
 there are infinitely many points $\beta \in E(\kbar)$ with 
$\hhat(\beta) < \varepsilon$ which are integral
with respect to $\alpha$.

\vskip .1 in

\begin{proof}
The argument is similar to the proof of Theorem \ref{Thm1}, 
but requires more machinery. We begin with some reductions.

First, after replacing $k$ by $k(\alpha)$, 
and $S$ by the set $S_{k(\alpha)}$ of places 
lying over $S$, we can assume that $\alpha \in k$.  

Second, after replacing $k$ by a finite extension $K/k$, 
and replacing $S$ with the set $S_K$ of places of $K$ lying above places in $S$, 
we can assume that $E$ has semi-stable reduction.  
Thus we can assume without loss of generality that for nonarchimedean $v$, 
either $E$  has good reduction, or $E$ is $k_v$-isomorphic to a Tate curve.  

Third, after enlarging $S$ if necessary, 
we can assume that $S$ contains all $v$ for which $|\Delta|_v \ne 1$.
In particular, we can assume that 
the model of $E$ defined by (\ref{FE1}) has good reduction 
for all $v \notin S$.

\vskip .1 in
We claim that if $\xi_n \in E(\kbar)_{\tors}$ is any torsion point, then 
\begin{equation} \label{FCOM6} 
\hhat(\alpha) \ = \ \frac{1}{[k(\xi_n):\QQ]} \sum_v 
\sum_{\sigma : k(\xi_n)/k \hookrightarrow \kbar_v} \lambda_v(\alpha-\sigma(\xi_n)) 
\ .   
\end{equation}
To see this, let $L$ be the Galois closure of $k(\xi_n)/k$. 
By (\ref{FCOM1}) and (\ref{FCOM3}), for each conjugate $\sigma(\xi_n)$, 
\begin{equation*} 
\hhat(\alpha) \ = \ \hhat(\alpha-\sigma(\xi_n)) 
\ = \ \frac{1}{[L:\QQ]} \sum_{\text{$w$ of $L$}} \lambda_w(\alpha - \sigma(\xi_n)) \ .
\end{equation*} 
Averaging over all embeddings $\sigma : L \hookrightarrow \kbar$, fixing 
an embedding $\kbar \hookrightarrow \kbar_v$ for each place $v$ of $K$, 
using (\ref{FCOM4}), 
and noting that there are only finitely many nonzero terms in each sum, we have 
\begin{eqnarray*} 
\hhat(\alpha) & = & \frac{1}{[L:k]} \sum_{\sigma : L/k \hookrightarrow \kbar} 
                      \frac{1}{[L:\QQ]} \sum_{\text{$w$ of $L$}}                  
                            \lambda_w(\alpha - \sigma(\xi_n)) \\
              & = & \frac{1}{[L:\QQ]} \sum_{\text{$v$ of $k$}} 
                    \sum_{\sigma : L/k \hookrightarrow \kbar_v}
                        \frac{1}{[L:k]} \sum_{w|v} [L_w:k_v] \cdot
                               \lambda_v(\alpha - \sigma(\xi_n)) \\
              & = & \frac{1}{[L:\QQ]} \sum_{\text{$v$ of $k$}} 
                    \sum_{\sigma : L/k \hookrightarrow \kbar_v} 
                               \lambda_v(\alpha - \sigma(\xi_n)) \ .
\end{eqnarray*}
Since each conjugate $\sigma(\xi_n)$ occurs $[L:k(\xi_n)]$ times in the final inner sum, this is equivalent to (\ref{FCOM6}).                                  

Suppose there were an infinite sequence of torsion points $\{\xi_n\}$ which were
$S$-integral with respect to $\alpha$.  

If $v \notin S$, our initial reductions assure that $E$ has good reduction at $v$.  By Proposition \ref{LocFormulas}.C and the integrality hypothesis, 
$\lambda_v(\alpha-\sigma(\xi_n)) = 0$ for each $n$ and $\sigma$.  It follows that
\begin{equation} \label{FCOM7} 
\hhat(\alpha) \ = \  \sum_{v \in S}  \frac{1}{[k(\xi_n):k]}
\sum_{\sigma : k(\xi_n)/k \hookrightarrow \kbar_v} \lambda_v(\alpha-\sigma(\xi_n)) 
\ .   
\end{equation} 

\vskip .1 in 
In the following two subsections, we will show that for each $v \in S$,  
\begin{equation} \label{FGF2}
\lim_{n \rightarrow \infty} 
\left( \frac{1}{[k(\xi_n):\QQ]} \sum_{\sigma : k(\xi_n)/k \hookrightarrow \kbar_v}
            \lambda_v(\alpha - \sigma(\xi_n)) \right) \ = \ 0 \ .
\end{equation} 
This will complete the proof of 
Theorem \ref{Thm2} for then, combining (\ref{FCOM7}) and (\ref{FGF2}) 
and letting $n \rightarrow \infty$ in (\ref{FCOM7}), 
we would have $\hhat(\alpha) = 0$,
contradicting the assumption that $\alpha$ is nontorsion.  

\subsubsection{\bf The Archimedean Case:}

\vskip .05 in

Let $v$ be an archimedean place of $k$.  
To simplify notation we view $k$ as embedded in $\CC$ and 
fix an isomorphism of $\kbar_v$ with $\CC$.
Thus, the way $k$ is embedded depends on the choice of $v$.  

To prove (\ref{FGF2}) we will need two facts:  David/Hirata-Kohno's theorem
on linear forms in elliptic logarithms, and a strong form of equidistribution
for torsion points.  

The following is a special case of (\cite{D-HK}, Theorem 1, p.31): 

\begin{proposition} \label{EllipticLF} {\rm (David/Hirata-Kohno)} 

Let $E/k$ be an elliptic curve defined over a number field
$k \subset \CC$.  Fix an isomorphism $\theta : \CC/\Lambda \cong E(\CC)$
for an appropriate lattice $\Lambda \subset \CC$.  Let $\omega_1, \omega_2$
be generators for $\Lambda$.  
Fix a non-torsion point $\alpha \in E(k)$
and let $a \in \CC$  be such that $\theta(a\text{\rm \ mod $\Lambda$}) = \alpha$. 

Then there is a constant $C = C(\alpha) > 0$ 
such that for all rational numbers $\ell_1/N, \ell_2/N$ 
with $\ell_1, \ell_2, N \in \ZZ$,  
\begin{equation*}
|a-(\frac{\ell_1}{N}\omega_1 + \frac{\ell_2}{N} \omega_2)| 
\ \ge \ e^{-C \max(1,\log(N))} \ .
\end{equation*} 
\end{proposition}

By Ullmo's theorem (\cite{Ull}), the Galois conjugates of the $\xi_n$ are 
equidistributed in $E(\CC)$.  As we will see, they are in fact strongly
equidistributed, in a sense analogous to that in Proposition \ref{Prop2}.  

If $\xi \in E(\kbar)_{\tors}$, write $\Gal(\kbar/k) \cdot \xi$
for the orbit $\{ \sigma(\xi) : \sigma \in \Gal(\kbar/k)\}$. 
For each set $U \subset E(\CC)$, write 
\begin{equation*}
N(\xi,U) \ = \ \#((\Gal(\kbar/k) \cdot \xi) \cap U) \ .
\end{equation*} 

Let $\cS \subset \CC$ be a bounded, convex, centrally symmetric set with $0$ 
in its interior.  For each $a \in \CC$ and $0 \le r \in \RR$, 
write $\cS(a,r) = \{ a + rz : z \in \cS \}$.  
For example, if $\cS = B(0,1)$ then $\cS(a,r) = B(a,r)$. 
 
Let $\Lambda \subset \CC$ be a lattice such that $E(\CC) \cong \CC/\Lambda$.   
Let $r_0 = r_0(\cS,\Lambda) > 0$ be the largest number such that 
$\cS(a,r)$  injects into $\CC/\Lambda \cong E(\CC)$ under the natural projection
for all $a \in \CC$ and all $0 \le r < r_0$.
Write $\cS_E(a,r)$ for the image of $\cS(a,r)$ in $E(\CC)$.  

\begin{proposition} \label{EquiProp} {\rm (Strong Equidistribution)} 
Let $k \subset \CC$ be a number field, and let $E/k$ be an elliptic curve.
Then the $\Gal(\kbar/k)$-conjugates of the torsion points in $E(\kbar)$
are strongly equidistributed in $E(\CC)$ in the following sense:

Let $\mu$ be the additive Haar measure on $E(\CC)$ with total mass $1$.  
Fix $\gamma$ with $0 < \gamma < 1/2$, 
and fix a bounded, convex, centrally symmetric set 
$\cS$ with $0$ in its interior.  
Then for each $r$ such that $\cS(a,r)$ injects into $E(\CC)$, 
and for all $\xi \in E(\kbar)_{\tors}$, 
\begin{equation*}
\frac{N(\xi,\cS_E(a,r))}{[k(\xi):k]}
 \ = \ \mu(\cS_E(a,r)) + O([k(\xi):k]^{-\gamma}) 
\end{equation*} 
where the implied constant depends only on $\cS$, $E$, and $\gamma$.
\end{proposition}

\vskip .1 in
The proof will be given in \S \ref{Sec2.2} below.  

\vskip .1 in
We can now complete the proof of (\ref{FGF2}) in the archimedean case.   
The argument is similar to the one in the proof of Theorem \ref{Thm1}.     
By Ullmo's theorem (\cite{Ull}), or by Proposition \ref{EquiProp} 
when $\cS$ has the shape of a period parallelogram 
(so $E$ can be tiled with sets $\cS_E(a,r)$), 
one knows that as $n \rightarrow \infty$ 
the discrete measures
\begin{equation*}
\mu_n \ = \ \frac{1}{[k(\xi_n):k]} 
    \sum_{\sigma: k(\xi_n)/k \hookrightarrow \CC} \delta_{\sigma(\xi_n)}(x) 
\end{equation*}
converge weakly to the Haar measure $\mu$ on $E(\CC)$  having total mass $1$.
Proving (\ref{FGF2}) is equivalent to showing that 
\begin{equation*}
\lim_{n \rightarrow \infty} \int_{E(\CC)} \lambda_v(\alpha-z ) \, d\mu_n(z) \ = \ 0  \ .
\end{equation*}

Choose a lattice $\Lambda \subset \CC$ such that $E(\CC) \cong \CC/\Lambda$,
and let $F$ be the area of a fundamental domain for $\Lambda$.  After scaling
$\Lambda$, if necessary, we can assume that $F = 1$.  After this normalization,
$\mu$ coincides with Lebesgue measure.  Let $\theta : \CC/\Lambda \cong E(\CC)$ 
be an isomorphism as in the David/Hirata-Kohno theorem, 
and let $a \in \CC$ be a point with $\theta(a\text{\rm \ mod $\Lambda$}) = \alpha$.  

Fix $\varepsilon > 0$ small enough that $B(a,\varepsilon)$ 
injects into $\CC/\Lambda$, 
and identify $B(a,\varepsilon)$ with its image 
$B_E(a,\varepsilon) = \theta(B(a,\varepsilon)) \subset E(\CC)$.
(In particular, identify $a$ with $\alpha$).  
Without loss, we can assume that $\varepsilon < 1/\pi$, 
so $\pi \varepsilon^2 < \varepsilon$.  
We will show that for all large $n$, 
\begin{equation} \label{FLZ1} 
|\int_{E(\CC)} \lambda_v(\alpha-z ) \, d\mu_n(z)| \ < \ 6 \varepsilon \ .
\end{equation} 

Put 
\begin{equation*}
\Theta_{\alpha,\varepsilon}(z) \ = 
\left\{ \begin{array}{cl}  
 -[ k_v : \RR ]  \log(|z-a|/\varepsilon) & \text{if $z \in B(a,r)$,} \\
 0                & \text{if $z \in E(\CC) \backslash B(a,r)$} 
\end{array} \right.  
\end{equation*}  
and note that
\begin{eqnarray*} 
0 & < & \int_{E(\CC)} \Theta_{\alpha,\varepsilon}(z) \, d\mu(z) \ = \ 
\int_{B(a,\varepsilon)} - [k_v: \RR] \log(|z-a|/\varepsilon) \, d\mu(z) \\
& = & [k_v : \RR ] \int_0^{\varepsilon} - 2 \pi t \log(t/\varepsilon) \, dt 
\ = \  [k_v : \RR] \frac{\pi \varepsilon^2}{2} \ < \ \varepsilon \ .
\end{eqnarray*}  
By Proposition \ref{LocFormulas}.A there is a continuous function $g_{\alpha,\varepsilon}(z)$ on $E(\CC)$ such that
\begin{equation*} 
\lambda_v(\alpha-z) \ = \ 
\Theta_{\alpha,\varepsilon}(z) + g_{\alpha,\varepsilon}(z) \ .
\end{equation*} 
Since $\int_{E(\CC)} \lambda_v(\alpha-z ) \, d\mu(z) = 0$  
(also by Proposition \ref{LocFormulas}.A), we get  
\begin{equation*}
|\int_{E(\CC)} g_{\alpha,\varepsilon}(z) \, d\mu(z)| \ = \ 
| \int_{B(a,\varepsilon)} -\Theta_{\alpha,\varepsilon}(z) \, d\mu(z)| 
\ < \ \varepsilon \ .
\end{equation*}
By weak convergence, it follows that for all sufficiently large $n$, 
\begin{equation} \label{FNN1} 
|\int_{E(\CC)} g_{\alpha,\varepsilon}(z) \, d\mu_n(z)| \ < \ 2 \varepsilon \ .
\end{equation} 
To complete the proof of (\ref{FLZ1}), it will suffice to show that for all sufficiently large $n$, 
\begin{equation} \label{FBAG}  
|\int_{B(a,r)} \log(|z-a|/\varepsilon) \, d\mu_n(z)| \ < \ 2 \varepsilon \ .
\end{equation} 

For this, put $D = D_n = \lceil [k(\xi_n):k]^{1/8} \rceil$,  
and subdivide $B(a,\varepsilon)$ into a 
disc $A_0(n) = B(a,\varepsilon/D)$ and annuli 
$A_{\ell}(n) = B(a,(\ell+1)\varepsilon/D) \backslash B(a,\ell \varepsilon/D)$ 
for $\ell =1, \ldots, D-1$.  

   For the central disc, we have 
$\mu(A_0(n)) = \pi \varepsilon^2/D^2 \le  \pi \varepsilon^2/[k(\xi_n):k]^{1/4}$.
Applying Proposition \ref{EquiProp} when $\cS$ is a disc, taking $\gamma = 3/8$,  
gives  
\begin{equation*}
N(\xi_n,A_1(n))/[k(\xi_n):k] \ \le \ 2 \mu(A_0(n)) 
\end{equation*} 
for all sufficiently large $n$.  If $\xi_n$ has order $N_n$,
the David/Hirata-Kohno theorem tells us that for each conjugate 
$\sigma(\xi_n) \in A_0(n)$ 
(where as before we are identifying $B(a,\varepsilon)$ with its image 
$\theta(B(a,\varepsilon)) \subset E(\CC)$) 
\begin{equation*}
|\log(|\sigma(\xi_n)-a|)| \ \le \ C \log(N_n) \ .
\end{equation*}
Using (\ref{FLL1}) and (\ref{FLL2}) below, one sees that   
$[k(\xi_n):k] \ge N_n^{1/2}$ for all sufficently large $n$.  Thus
 $0 \le |\log(|\sigma(\xi_n)-\alpha|)| \ \le \ 2C \log([k(\xi_n):k])$
and 
\begin{equation} \label{FQQ1} 
0 \ \le \ |\int_{A_0(n)} \log(|z-\alpha|) \, d\mu_n(z)| 
  \ \le \  4\pi \varepsilon^2 C \cdot
  \frac{\log([k(\xi_n):k])}{[k(\xi_n):k]^{1/4}} \ < \ \varepsilon
\end{equation}
for all sufficiently large $n$.  

For each annulus $A_{\ell}(n)$, $\ell = 1, \ldots, D-1$, one has 
\begin{equation*} 
\mu(A_{\ell}(n)) \ = \ \pi (2 \ell + 1) \varepsilon^2 /D^2 
\ \cong \ \pi (2 \ell+1) \varepsilon^2 /[k(\xi):k]^{1/4} \ .
\end{equation*}   
Since $A_{\ell}(n)$ is the difference of two sets to which 
Proposition \ref{EquiProp} applies, 
we find as above that for sufficiently large $n$,   
\begin{equation*}
N(\xi_n,A_{\ell}(n))/[k(\xi_n):k] \ \le \ 2 \mu(A_{\ell}(n)) \ .
\end{equation*} 
Note that on $A_{\ell}(n)$,  
$|\log(|z-\alpha|/\varepsilon)| \ \le \ -\log(\ell/D)$.  
Summing over these annuli, and bounding the resulting Riemann sum
by an integral, we find that
\begin{eqnarray}
|\int_{B(a,\varepsilon) \backslash A_0(n)} \log(|z-a|/\varepsilon)\, d\mu_n(z)| 
  & \le & \sum_{\ell=1}^{D-1} -\log((\frac{\ell \varepsilon}{D})/\varepsilon) 
                 \cdot 2 \mu(A_{\ell}(n))   \notag \\ 
  & < & 2 \cdot \int_{B(a,\varepsilon)} -2 \pi t \log(t/\varepsilon) \, dt \notag \\
  & = & \pi \varepsilon^2 \ < \ \varepsilon \ . \label{FNQ1}
\end{eqnarray}  
Combining   (\ref{FQQ1}) and (\ref{FNQ1}) gives (\ref{FBAG}), which  
completes the proof of (\ref{FGF2}) in the archimedean case (assuming
Proposition~\ref{EquiProp}).

\subsubsection{\bf The Nonarchimedean Case:}
\vskip .05 in

In the nonarchimedean case, the proof of (\ref{FGF2}) depends on 
a well-known result of Cassels on the denominators
of torsion points (see \cite{Sil1}, Theorem 3.4, p.177).  
Write $\bcO_v$ for the ring of integers of $\kbar_v$.  

\begin{proposition} \label{CasselsThm} {\rm (Cassels)} 
\nopagebreak

Let $k_v$ be a local field of characteristic $0$ 
and residue characteristic $p  > 0$, and let $E/k_v$
be an elliptic curve defined by a Weierstrass equation 
\begin{equation*}
y^2 + a_1 xy + a_3 y \ = \ x^3 + a_2 x^2 + a_4 x + a_6
\end{equation*}
whose coefficients belong to $\cO_v$ 
$($N. B. the Weierstrass equation need not be minimal.$)$  
Let $P \in E(\kbar_v)_{\tors}$ be a point of exact order $m \ge 2$.

$($A$)$  If $m$ is not a power of $p$, then $x(P), y(P) \in \overline{\cO}_v$. 

$($B$)$  If $m = p^n$, then $x(P) = a/D^2$, $y(P) = b/D^3$ 
where $a, b, D \in \overline{\cO}_v$ and 
\begin{equation*}
\ord_v(D) \ \le \  \frac{\ord_v(p)}{p^n - p^{n-1}}  \ .
\end{equation*}
\end{proposition} 

\begin{proof} 
Silverman (\cite{Sil1}, Theorem 3.4) 
states the theorem for torsion points belonging to 
$E(k_v)$, with $a, b, D \in k_v$ in part B) and $D$ satisfying
\begin{equation} \label{FVD1} 
\ord_v(D) \ = \ \big\lfloor \frac{\ord_v(p)}{p^n - p^{n-1}} \big\rfloor \ .
\end{equation}
Since the Weierstrass equation for $E$ need not be minimal, 
we can replace $k_v$ by an arbitrary finite extension $L_w/k_v$, 
and if $e_{w/v}$ is the ramification index of $L_w/k_v$, 
then for $P \in E(L_w)_{\tors}$ and $a, b, D \in L_w$,
(\ref{FVD1}) becomes 
\begin{equation} \label{FVD2} 
\ord_v(D) \ = \ \frac{1}{e_{w/v}} 
\cdot \big\lfloor \frac{e_{w/v} \ord_v(p)}{p^n - p^{n-1}} \big\rfloor \ .
\end{equation} 
This yields the result for all $P \in E(\kbar_v)_{\tors}$.  
\end{proof} 

\vskip .1 in
\begin{corollary} \label{UpperBd} 
Let $E/k_v$ be an elliptic curve defined over a nonarchimedean local field.
Then for each nontorsion point $\alpha \in E(\kbar_v)$:

$($A$)$  There is a number $M$ such that for all $\xi \in E(\kbar_v)_{\tors}$,
\begin{equation*}
\lambda_v(\alpha-\xi) \ \le \ M \ .
\end{equation*}

$($B$)$  If $E$ has good reduction, then for each $\varepsilon > 0$, there
are only finitely many $\xi \in E(\kbar_v)_{\tors}$ with 
$\lambda_v(\alpha-\xi) \ > \ \varepsilon $. 
 
If $E$ is a Tate curve, then for each $\varepsilon > 0$, there
are only finitely many $\xi \in E(\kbar_v)_{\tors}$ with 
$\lambda_v(\alpha-\xi) > \varepsilon + \frac{1}{12} (-\log(|\Delta(E)|_v))$.  
\end{corollary}

\begin{proof}
After a finite base extension, we can assume that $E$ either has good reduction
or is a Tate curve.  Since (B) implies (A), it suffices to prove (B).  
Fix $\varepsilon > 0$.  

First suppose $E$ has good reduction.  
Then $\lambda_v(x-y) = -\log(\|x,y\|_v)$, where $\|x,y\|_v$
is the spherical distance on the minimal Weierstrass model for $E/k_v$.
If $\xi_1, \xi_2 \in E(\kbar_v)_{\tors}$
satisfy $\lambda_v(\alpha-\xi_i) > \varepsilon$, then 
$\|\xi_1,\alpha\|_v$, $\|\xi_2,\alpha\|_v < (Nv)^{-\varepsilon}$
where $Nv$ is the order of the residue field of $\cO_v$. 
By the the ultrametric inequality for the spherical distance
(\cite{R1}, \S1.1), $\|\xi_1,\xi_2\|_v < (Nv)^{-\varepsilon}$.  
By translation invariance  
$\|\xi_1-\xi_2,0\|_v  < (Nv)^{-\varepsilon}$.  
Put $\xi := \xi_1-\xi_2$.  By the definition of the spherical distance, 
if $x, y$ are the coordinate functions in the minimal Weierstrass model,
\begin{equation*}
-\log(\|\xi,0\|_v) = \min(\ord_v(x(\xi)),\ord_v(y(\xi))) \cdot \log(Nv)  \ .
\end{equation*}
By Cassels' theorem, there only finitely many torsion points for which
\begin{equation*}
\min(\ord_v(x(\xi)),\ord_v(y(\xi))) \ > \ \varepsilon/\log(Nv) \ . 
\end{equation*}  

Next suppose $E$ is a Tate curve.  
Fix a Tate isomorphism 
$E(\kbar_v) \cong \kbar_v^{\times}/q^{\ZZ}$ where $|q|_v = |\Delta(E)|_v < 1$,
and let $y^2 +xy = x^3 + a_4(q)x + a_6(q)$ 
be the corresponding Weierstrass equation.  
Let $a, u_1, u_2 \in \kbar_v^{\times}$  correspond to $\alpha, \xi_1, \xi_2$
respectively;  we can assume  that 
$|q|_v < |a|_v, |u_1|_v, |u_2|_v \le 1$.  By the formula for 
$\lambda_v(x-y)$ in Proposition \ref{LocFormulas}.B, 
if $\lambda_v(\alpha-\xi_i) > \varepsilon + \frac{1}{12} (-\log(|\Delta(E)|_v))$,
then $|\, a|_v = |u_1|_v = |u_2|_v$ and 
\begin{equation*}
-\log(|1-a^{-1}u_i|_v) \ = \ 
\ord_v(1-a^{-1}u_i) \cdot \log(Nv) \ > \ \varepsilon  \ .
\end{equation*} 
Put $\xi = \xi_1-\xi_2$
and $u = u_2^{-1}u_1$.  Then $\xi$ corresponds to $u$ under the Tate 
isomorphism, and $\ord_v(1-u) > \varepsilon/\log(Nv)$.  
By the formulas for $x(\xi)$, $y(\xi)$ in (\cite{Sil2}, p. 425), 
$\ord_v(x(\xi)) = 2 \, \ord_v(1-u) $ and $\ord_v(y(\xi)) = 3 \, \ord_v(1-u)$.   
Again by Cassels' theorem, only finitely many torsion points $\xi$ 
can satisfy $\min(\ord_v(x(\xi)),\ord_v(y(\xi))) > \varepsilon/\log(Nv)$.  
\end{proof} 

\vskip .1 in
We can now prove (\ref{FGF2}) when $E$ has good reduction at $v$.  

Fix $\varepsilon > 0$.  
Let $M$ be the upper bound in Corollary \ref{UpperBd}.A,
and let $N$ be the number of points
$\xi \in E(\kbar_v)_{\tors}$ with $\lambda_v(\alpha-\xi) > \varepsilon$ 
given by Corollary \ref{UpperBd}.B.  
For all sufficiently large $n$, $MN/[k(\xi_n):k] < \varepsilon$, giving 
\begin{eqnarray*}
0 & \le & \frac{1}{[k(\xi_n):k]} 
         \sum_{\sigma : \kbar/k \hookrightarrow \kbar_v} 
            \lambda_v(\alpha-\sigma(\xi_n)) \\
  & \le & \frac{([k(\xi_n):k] - N)}{[k(\xi_n):k]} \cdot \varepsilon 
            +  \frac{N}{[k(\xi_n):k]} \cdot M \ < \ 2 \varepsilon \ .
\end{eqnarray*}
Thus 
\begin{equation*}
\lim_{n \rightarrow \infty}  \frac{1}{[k(\xi_n):k]} 
          \sum_{\sigma : \kbar/k \hookrightarrow \kbar_v} 
            \lambda_v(\sigma(\xi_n)-\alpha)  \ = \ 0 \ . 
\end{equation*}   

\vskip .1 in  

To prove (\ref{FGF2}) when $E$ is a Tate curve at $v$, 
we will need the following
equidistribution theorem of Chambert-Loir (\cite{CL}, Corollaire 5.5). 

Fix a Tate isomorphism $E(\kbar_v) \cong \kbar_v/q^{\ZZ}$, 
put $L = \ZZ \cdot \ord_v(q) \subset \RR$, and  
define a ``reduction map'' 
$r : E(\kbar) \rightarrow \RR/L$ by setting
$r(P) = \ord_v(a) \text{\ (mod $L$) }$ 
if $P \in E(\kbar_v)$ corresponds to $a \in \kbar_v^{\times}$.  

For each global point $P \in E(\kbar)$, 
define a measure $\mu_{P,v}$ on $\RR/ L$ by
\begin{equation*}
\mu_{P,v}(z)  \ = \ \frac{1}{[k(P):k]} 
      \sum_{\sigma: \kbar/k \hookrightarrow \kbar_v} 
     \delta_{r(\sigma(P))}(z)  
\end{equation*}
and let $\mu_v$ be the Haar measure on $\RR/ L$ with
total mass $1$.     

\begin{proposition} \label{Chambert} {\rm (Chambert-Loir)}

For each sequence of points $\{P_n\}$ in $E(\kbar)$ 
with $\hhat(P_n) \rightarrow 0$, 
the  sequence of measures $\{\mu_{P_n,v}\}$
converges weakly to $\mu_v$.  
\end{proposition} 

\vskip .1 in 
We can now prove (\ref{FGF2}) when $E$ is a Tate curve.  Let $\{\xi_n\}$
be a sequence of torsion points which are $S$-integral 
with respect to $\alpha$.  

Fix $\varepsilon > 0$.  
Let $M$ be the upper bound in Corollary \ref{UpperBd}.A.  
Put $a = r(\alpha)$ and
let $\delta > 0$ be such that $\mu((a-\delta,a+\delta)) < \varepsilon/M$, 
where by abuse of notation we identify a sufficiently short interval
in $\RR$ with its image in $\RR/L$.    
By Chambert-Loir's theorem, 
$\mu_{\xi_n,v}((a-\delta,a+\delta)) < 2 \varepsilon/M$
for all sufficiently large $n$.

By the formulas in Proposition \ref{LocFormulas}.B, 
$\int_{\RR/L} \tlambda_v(z) \, d\mu_v(z) = 0$
and 
\begin{eqnarray*}
& & \big| \frac{1}{[k(\xi_n):k]} 
       \sum_{\sigma : \kbar/k \hookrightarrow \kbar_v} 
            \lambda_v(\sigma(\xi_n)-\alpha) \big| \\
& & \qquad  \ \le \  |\int_{\RR/L} 
           \tlambda_v(z-a) \, d\mu_{\xi_n,v}(z) | 
\ + \ M \cdot \mu_{\xi_n,v}((a-\delta,a+\delta)) \ .
\end{eqnarray*}
For sufficiently large $n$ the right side is at most $3 \varepsilon$.  
Hence        
\begin{equation*}
\lim_{n \rightarrow \infty}  \frac{1}{[k(\xi_n):k]} 
          \sum_{\sigma : \kbar/k \hookrightarrow \kbar_v} 
            \lambda_v(\sigma(\xi_n)-\alpha) \ = \ 0 \ . 
\end{equation*}  
This completes the proof of Theorem \ref{Thm2}. 
\end{proof} 
 
\vskip .1 in
\vskip .05 in
Several results in the literature use methods similar to ours, though 
none of them yields Theorem \ref{Thm2}:    
 
J.~Cheon and S.~Hahn \cite{Ch-Han} proved an elliptic curve analogue 
of Schinzel's theorem \cite{Schinzel}.  
Likewise, Everest and B. N\'i Flath\'uin \cite{EnF} 
evaluate `elliptic Mahler measures' 
in terms of limits involving division polynomials,
obtaining results similar to (\ref{BB2}).     
They use David/Hirata-Kohno's theorem on elliptic logarithms in place of Baker's theorem, much as we do.  
     
More recently, L. Szpiro and T. Tucker \cite{SZ-T} proved that 
local canonical heights for a dynamical system can be evaluated 
by taking limits over `division polynomials' for the dynamical system.  
(These polynomials have periodic points as their roots).  
Their work uses Roth's theorem rather than Baker's or David/Hirata-Kohno's theorem.  
It would be interesting to see if this could be brought to bear
on Conjecture \ref{Conj2}. 

\vskip .1 in
\subsection{\bf Strong equidistribution for torsion points on elliptic curves.} 
\label{Sec2.2} 

\vskip .1 in
We will now prove Proposition~\ref{EquiProp}, the strong equidistribution theorem for torsion points 
on elliptic curves which was used in the proof of Theorem~\ref{Thm2}.

\vskip .1 in 
\begin{proof} of Proposition \ref{EquiProp}.  

The proof breaks into two cases, depending on whether $E$ has complex
multiplication or not.  
\vskip .1 in

First suppose $E$ does not have complex multiplication. 
 
As usual, the action of $\Gal(\kbar/k)$ on $E(\kbar)_{\tors}$ gives 
a homomorphism
\begin{equation*}
\eta : \Gal(\kbar/k) \rightarrow \lim_{\longleftarrow} \GL_2(\ZZ/N\ZZ) 
\cong \prod_p \GL_2(\ZZ_p) \ .
\end{equation*}
By Serre's theorem (\cite{Serre}, Th\'eor\`eme 3), 
the image of $\Gal(\kbar/k)$ in $\prod_p \GL_2(\ZZ_p)$ is open. Thus
there is a number $Q$ such that $\Im(\eta)$ contains
the subgroup 
\begin{equation*}
\prod_{p|Q} (1+Q M_2(\ZZ_p)) \times 
\prod_{p \!\not \, \mid Q} \GL_2(\ZZ_p) \ .
\end{equation*}
Let $G_Q \subset \Gal(\kbar/k)$ be the pre-image of this subgroup.  

Let $\xi \in E(\kbar)_{\tors}$ have order $N$, and put $Q_N = \gcd(Q,N)$.
For suitable right coset representatives 
$\sigma_1, \ldots, \sigma_T$ of $G_Q$ in $\Gal(\kbar/k)$, the Galois
orbit $\Gal(\kbar/k) \cdot \xi$ decomposes as a 
disjoint union of $G_Q$-orbits: 
\begin{equation*}
\Gal(\kbar/k) \cdot \xi 
\ = \ \bigcup_{i=1}^T G_Q \cdot \sigma_i(\xi) \ .
\end{equation*}
Since $G_Q$ is normal in $\Gal(\kbar/k)$,
the orbits $G_Q \cdot \sigma_i(\xi) = \sigma_i(G_Q \cdot \xi)$ 
all have the same size.   
Thus $[k(\xi):k] = T \cdot \#(G_Q \cdot \xi)$.  By considering the 
action of $G_Q$ on the $p$-parts of $\xi$,
one sees that
\begin{eqnarray}
\#(G_Q \cdot \xi) & = & \prod_{p |Q_N} p^{2(\ord_p(N)-\ord_p(Q_N))} 
            \cdot \prod_{\substack{ p|N \\ p \!\not \, \mid Q_N }}
                     p^{2\ord_p(N)} (1 - \frac{1}{p^2}) \notag \\
            & = & \frac{N^2}{Q_N^2} \cdot 
            \prod_{\substack{ p|N \\ p \!\not \, \mid Q }}
                     (1 - \frac{1}{p^2}) \ . \label{FDV1} 
\end{eqnarray} 
Indeed, let $\xi_p$ be the $p$-component of $\xi$ in 
$E[N] \cong \prod_{p|N} (\ZZ/p^{\ord_p(N)}\ZZ)^2$.  
Identify $\xi_p$ with an element of $(\ZZ/p^{\ord_p(N)}\ZZ)^2$,
and note that it is a generator for that group. 
If $p|Q_N$, the image of $G_Q$ in $\GL_2(\ZZ/p^{\ord_p(N)}\ZZ)$ 
is $I + p^{\ord_p(Q_N)}M_2(\ZZ/p^{\ord_p(N)}\ZZ)$, and
\begin{equation*}
G_Q \cdot \xi_p \ =  
            \xi_p + p^{\ord_p(Q_N)} \cdot (\ZZ/p^{\ord_p(N)}\ZZ)^2  \ .
\end{equation*}
On the other hand, if $p \!\not \, \mid Q_N$, the image of $G_Q$
in $\GL_2(\ZZ/p^{\ord_p(N)}\ZZ)$ is the full group, so
\begin{equation*}                        
G_Q \cdot \xi_p \ = \ (\ZZ/p^{\ord_p(N)}\ZZ)^2 \backslash  
                         p \cdot (\ZZ/p^{\ord_p(N)}\ZZ)^2 \ .
\end{equation*}                   

Write $\Lambda_N = \frac{1}{N} \Lambda$, fix $\sigma_i$,  
and let $x \in \Lambda_N$ correspond to $\sigma_i(\xi)$.
Since $E[N] \cong \Lambda_N/\Lambda$, the considerations above show   
there is a one-to-one correspondence 
between elements of  $G_Q \cdot \sigma_i(\xi)$,
and cosets $y + \Lambda$ for $y \in \Lambda_N$ such that 
$y-x \in Q_N\Lambda_N$ and $y + \Lambda$ has exact order
$N$ in $\Lambda_N/\Lambda$.  Equivalently, 
$y-x \in Q_N \Lambda_N$ and $y \notin p\Lambda_N$ 
for each prime $p$ dividing $N$ but not $Q$.  

Let $p_1, \ldots, p_R$ be the primes dividing $N$ but not $Q$;  
if there are no such primes, take $p_1 \cdots p_R = 1$.
Since $Q_N$ and $p_1, \cdots, p_R$ are pairwise coprime, there is an 
$x_0 \in \Lambda_N$ such that 
$x_0 \equiv x \text{\ (mod $Q_N \Lambda_N$)}$ 
and $x_0 \equiv 0 \text{\ (mod $p_1 \cdots p_R \Lambda_N$)}$. 
Then $y-x \in Q_N \Lambda_N$ if and only 
if $y \in x_0 + Q_N \Lambda_N$, and $y \in p_i \Lambda_N$
if and only if $y \in x_0 + p_i \Lambda_N$.  
Note that if $D|p_1 \cdots p_R$ then 
$Q_N\Lambda_N \cap D\Lambda_N = Q_N D \Lambda_N$.
Take $a \in \CC$ and $0 < r \le r_0$. 
Using the fact that $\cS(a,r)$ injects into $\CC/\Lambda$ 
and applying inclusion-exclusion, we obtain
\begin{eqnarray}
& & \#(G_Q \cdot \sigma_i(\xi) \cap \cS_E(a,r)) \label{FAC1} \\
& & \qquad \qquad \ = \ 
\sum_{D|p_1 \cdots p_R} (-1)^{\lambda(D)} \cdot 
              \# (  (x_0 +  Q_N D\Lambda_N) \cap \cS(a,r) ) \notag
\end{eqnarray}      
where $\lambda(D)$ is number of primes dividing $D$.  

Let $\cF$ be a fundamental domain for $\Lambda$;
we can assume $\cF$ is bounded and contains $0$.  
Let $C$ be such that $\cF \subset \cS(0,C)$.   
Note that since $\cS$ is convex,
if $z_1 \in \cS(a_1,r_1)$ and $z_2 \in \cS(a_2,r_2)$,  
then $z_1+z_2 \in \cS(a_1+a_2,r_1+r_2)$.  
Put $F = \area(\cF)$, $S = \area(\cS)$;
then $\area(t\cF) = t^2 F$ and $\area(\cS(a,r)) = r^2 S$.

Each lattice $Q_N D \Lambda_N$ is homothetic to $\Lambda_N$, and hence 
has fundamental domain $(Q_ND/N) \cdot \cF \subset \cS(0,C \cdot Q_ND/N)$.
Write $t = Q_ND/N$, so $Q_N D \Lambda_N = t\Lambda$  
and $t\cF \subset \cS(0,tC)$.  As $y$ runs over $x_0+t\Lambda$, 
the sets $y+t\cF$ are pairwise disjoint and cover $\CC$.
If $y \in \cS(a,r)$, then $y+t\cF \subset \cS(a,r+tC)$.   Hence  
\begin{eqnarray}
\# ( (x_0 + t\Lambda) \cap \cS(a,r) )  
     & \le & \frac{\area(\cS(a,r+tC))}{\area(t \cF)} \label{FBA1} \\
     & = & \frac{ r^2 S}{F} \cdot \frac{1}{t^2}  
         + \frac{2CSr}{F} \cdot \frac{1}{t} + \frac{C^2 S}{F} \ . \notag 
\end{eqnarray}         
Similarly, if $r > tC$, take $z \in \cS(a,r-tC)$, 
and let $y \in x_0+t \Lambda$ be such that $z \in y + t \cF$.  Then
$z-y \in t \cF$, so $z-y \in \cS(0,tC)$, and since $\cS$ is centrally
symmetric $y-z \in \cS(0,tC)$.  Thus $y = z + (y-z) \in \cS(a,r)$.  
It follows that $\cS(a,r-tC) \subset 
\bigcup_{y \in (x_0 + t\Lambda) \cap \cS(a,r) } (y +t\cF)$, so 
\begin{eqnarray} 
\# ( (x_0 + t\Lambda) \cap \cS(a,r) )   
     & \ge & \frac{\area(\cS(a,r-tC))}{\area(t \cF)}  \label{FBA2} \\
     & > & \frac{ r^2 S}{F} \cdot \frac{1}{t^2}  
         - \frac{2CSr}{F} \cdot \frac{1}{t} - \frac{C^2 S}{F} \ . \notag                     
\end{eqnarray}
If $r \le tC$ then the right side of (\ref{FBA2}) is negative, so the 
inequality between the first and last quantities holds trivially.  

Replacing $t$ by its value $Q_N D/N$ and combining (\ref{FBA1}), (\ref{FBA2}),
we obtain 
\begin{eqnarray}
& & \big|\# ((x_0 + Q_N D\Lambda_N) \cap \cS(a,r)) 
 - \frac{\area(\cS(a,r))}{\area(\cF)} \cdot \frac{N^2}{Q_N^2 D^2} \big| \notag \\
& & \qquad \qquad \qquad \qquad  
       \ \le \  \frac{2  C S r}{F} \cdot \frac{N}{Q_N D} 
                + \frac{ C^2S }{F} \ .  \label{FBA3}
\end{eqnarray}  

Inserting (\ref{FBA3}) in the inclusion-exclusion relation (\ref{FAC1})
and summing over all $\sigma_i(\xi)$, $i = 1, \ldots, T$, we find 
\begin{eqnarray*}
N(\xi,\cS_E(a,r)) 
& = & \frac{\area(\cS(a,r))}{\area(\cF)} \cdot \frac{T N^2}{Q_N^2} 
   \prod_{ p|N , p \!\not \, \mid Q } (1-\frac{1}{p^2}) \\
& & \ + \ \theta\big(\frac{2 C S r}{F} \cdot \frac{TN}{Q_N}
 \prod_{ p|N , p \!\not \, \mid Q } (1+\frac{1}{p})\big)
    + \ \theta\big( \frac{C^2 S}{F} \cdot T 2^R \big)
\end{eqnarray*}
where as before, 
$\theta(x)$ denotes a quantity with $-x \le \theta(x) \le x$.
By (\ref{FDV1}), 
\begin{equation} \label{FGJ1} 
[k(\xi):k] \ = \ T \cdot \#(G_Q \cdot \xi) \ = \ 
 \frac{TN^2}{Q_N^2} \prod_{p|N, p \!\not \, \mid Q} (1-1/p^2) \ .  
\end{equation}
Since $r \le r_0$, it follows that 
\begin{eqnarray*}
\frac{N(\xi,\cS_E(a,r))}{[k(\xi):k]} & = & \frac{\area(\cS(a,r))}{\area(\cF)} 
 \ + \ \theta\big(\frac{2 C S r_0}{F} \cdot 
    \frac{Q_N}{N \prod_{ p|N , p \!\not \, \mid Q } (1-\frac{1}{p})} \big)\\
& & \qquad \qquad + \ \theta\big(  \frac{ C^2 S}{F} \cdot \frac{2^R Q_N^2}
      {N^2 \prod_{ p|N , p \!\not \, \mid Q } (1-\frac{1}{p^2})} \big) \ .
\end{eqnarray*}

Here $ \area(\cS(a,r))/\area(\cF) = \mu(\cS_E(a,r))$.   
Note that $T \le \#(\GL_2(\ZZ/Q\ZZ))$ is bounded, $Q_N \le Q$ is bounded, 
and $N \prod_{p|N}(1-\frac{1}{p}) \ge N^{1-\varepsilon}$ 
for each $\varepsilon > 0$ and each sufficiently large $N$.
Using (\ref{FGJ1}) and the fact that 
$1 \ge \prod_{ p|N , p \!\not \, \mid Q } (1-\frac{1}{p^2}) \ge 1/\zeta(2)$
one sees that the first error term is $O_{\gamma}([k(\xi):k]^{-\gamma})$
for each $\gamma < 1/2$.  Similarly, $2^R \le d(N) \le N^{\varepsilon}$
for each $\varepsilon > 0$ and each sufficiently large $N$. 
Thus the second error term is negligible in comparison with the first.

This completes the proof of Proposition \ref{EquiProp}
when $E$ does not have complex multiplication. 

\vskip .1 in
Now suppose $E$ has complex multiplication.  Let $K$ be the $CM$ field,
and let $\cO \subset \cO_K$ be the order corresponding to $E$.  After
enlarging $k$ if necessary, we can assume that $K \subset k$.  
Let $\Lambda \subset \CC$ be a lattice such that $E \cong \CC/\Lambda$.
Without loss of generality, we can assume that $\Lambda \subset K$.  
Fix an analytic isomorphism $\vartheta : \CC/\Lambda \cong E(\CC)$.  

By the theory of complex multiplication 
(see \cite{Sh}, \cite{Lang1}, or \cite{Sil2}, Chapter II), 
$E(\kbar)_{\tors}$ is rational over $k^{ab}$, the maximal abelian extension of $k$.  
Let $k_{\AA}^{\times}$ be the idele ring of $k$, 
and for $s \in k_{\AA}^{\times}$ let $[s,k]$ be the Artin map  
acting on $k^{ab}$.  Given $\sigma \in \Gal(\kbar/k)$, 
take $s \in k_{\AA}^{\times}$ with $\sigma\vert_{k^{ab}} = [s,k]$,
and put $w = N_{k/K}(s) \in K_{\AA}^{\times}$.   There is an action
of $K_{\AA}^{\times}$ on lattices, defined semi-locally, which 
associates to $w$ and $\Lambda$ a new lattice $w^{-1}\Lambda$. 
This action extends to a map $w^{-1}: K/\Lambda \rightarrow K/w^{-1} \Lambda$.
There is also a homomorphism $\psi: k_{\AA}^{\times} \rightarrow K^{\times}$,
the `gr\"ossencharacter' of $E$, which has the property that 
$\psi(s) N_{k/K}(s)^{-1} \Lambda = \Lambda$.  
Put $\kappa = \psi(s) \in K^{\times}$.

With this notation, there is a commutative diagram: 
\begin{equation*}
\begin{array}{rcccccl}
       & K/\Lambda & \hookrightarrow & \CC/\Lambda 
             & \stackrel{\vartheta}{\longrightarrow} & E(\kbar)_{\tors} & \\
        w^{-1}  & \downarrow &  &  & & \downarrow & \sigma \\          
      & K/w^{-1} \Lambda & \hookrightarrow & \CC/w^{-1}\Lambda 
              & \longrightarrow  & E(\kbar)_{\tors} & \\
    \kappa  & \downarrow &  &  & & \downarrow & id \\ 
    & K/\Lambda & \hookrightarrow & \CC/\Lambda 
            & \stackrel{\vartheta}{\longrightarrow} & E(\kbar)_{\tors} & 
                \end{array}
\end{equation*}
in which the vertical arrows on the left are multiplication by $w^{-1}$
and $\kappa$ respectively, and those on the right are the Galois action
(see \cite{Sh}, Proposition 7.40, p.211, or \cite{Lang1}, Theorem 8, p.137).  
Note that the same analytic isomorphism $\vartheta$ 
appears in the top and bottom rows.
Thus, if $\xi \in E(\kbar)_{\tors}$ corresponds to $x \in K/\Lambda$, 
and $\sigma|_{k^{ab}} = [s,k]$, then  
\begin{equation*}
\sigma(\xi) \ = \ \vartheta(\psi(s) N_{k/K}(s)^{-1} x) \ .
\end{equation*} 
This gives an explicit description of the Galois action on torsion points  
in terms of adelic ``multiplication''.  

The action of $K_{\AA}^{\times}$ in the diagram is as follows.  
Let $L \subset K$ be a lattice.  For each rational prime $p$ of $\QQ$, 
write $L_p = L \otimes_{\ZZ} \ZZ_p$ and $K_p = K \otimes_{\QQ} \QQ_p$;  
if $w \in K_{\AA}^{\times}$, let $w_p$ be its $p$-component.  
Then $w_p^{-1}L_p$ is a $\ZZ_p$-lattice in $K_p$.  There is a unique 
lattice $M \subset K$ such that $M_p = w_p^{-1}L_p$ for each $p$ 
(\cite{Lang1}, Theorem 8, p.97), and $w^{-1}L$ is defined to be $M$.
Likewise, if $x \in K/L$, lift it to an element of $K \subset K_{\AA}$ and  
write $x_p \in K_p$ for its $p$-component;  there is a $y \in K$ such 
that $w_p^{-1} x_p \text{\ (mod $w^{-1}L_p$)} = y \text{\ (mod $M_p$)}$
for each $p$, and  $w^{-1}(x \text{\ (mod $L$)})$ is defined to be 
$y \text{\ (mod $M$)}$.

The order $\cO$ has the form $\cO = \ZZ + c \cO_K$ 
for some integer $c \ge 1$, and $c$ is called the conductor of $\cO$.   
The lattice $\Lambda$ is a proper $\cO$-lattice, meaning that 
$\cO = \{x \in K : x \Lambda \subset \Lambda\}$.  For any order $\cO$, 
there are only finitely many homothety classes of proper $\cO$-lattices 
(\cite{Lang1}, Theorem 7, p.95). 
Write $\cO_p = \cO \otimes_{\ZZ} \ZZ_p$ and 
$\cO_{K,p} = \cO_K \otimes_{\ZZ} \ZZ_p$.  If $p \!\not \, \mid c$,
then $\cO_p = \cO_{K,p} \cong \prod_{\fp |p} \cO_{K,\fp}$, 
where $\fp$ runs over the primes of $K$ lying over $p$, and $\cO_{K,\fp}$
is the completion of $\cO_K$ at $\fp$.   

The kernel $U$ of the gr\"ossencharacter 
$\psi : k_{\AA}^{\times} \rightarrow K^{\times}$ is open in $k_{\AA}^{\times}$, 
so its image $W = N_{k/K}(U) \subset K_{\AA}^{\times}$ is open.  Thus there 
is an integer $Q \ge 1$ such that for each $p|Q$, the subgroup 
$1 + Q \cO_{K,p} \subset \cO_{K,p}^{\times}$ is contained in $W_p$ 
and for each $p \!\not \, \mid Q$, $\cO_{K,p}^{\times} \subset W_p$.  
If $w \in W$, then $w^{-1} \Lambda = \Lambda$, so $w_p \in \cO_p^{\times}$.
Hence $c|Q$.   

Noting that $\cO_p = \cO_{K,p}$ if $p \!\not \, \mid Q$,
let $W_Q \subset \ K_{\AA}^{\times}$ be the subgroup 
\begin{equation*} 
\CC^{\times} \times \prod_{p |Q} (1 + Q \cO_p) \times 
               \prod_{p \!\not \, \mid Q} \cO_p^{\times}
       \ \subset \  W \ ,
\end{equation*}
and let $U_Q$ be its preimage in $k_{\AA}^{\times}$ under the
norm map. Put  
\begin{equation*}
G_Q \ = \ \{ \sigma \in \Gal(\kbar/k) : \sigma|k^{ab} = [s,k] 
                         \text{\ for some $s \in U_{Q}$} \} \ .
\end{equation*}                        
Then $G_Q$ is open and normal in $\Gal(\kbar/k)$.
  
Fix $\xi \in E(\kbar)_{\tors}$.  Suppose $\xi$ has order $N$; 
put $Q_N = \gcd(Q,N)$.
For suitable right coset representatives 
$\sigma_1, \ldots, \sigma_T$ of $G_Q$ in $\Gal(\kbar/k)$, the 
orbit $\Gal(\kbar/k) \cdot \xi$ decomposes as a 
disjoint union of $G_Q$-orbits: 
\begin{equation*}
\Gal(\kbar/k) \cdot \xi 
\ = \ \bigcup_{i=1}^T G_Q \cdot \sigma_i(\xi) \ .
\end{equation*}
As before, the orbits $G_Q \cdot \sigma_i(\xi) = \sigma_i(G_Q \cdot \xi)$ 
all have the same size, and  $[k(\xi):k] = T \cdot \#(G_Q \cdot \xi)$.  

Let $\xi$ correspond to $x + \Lambda \in K/\Lambda$.  
Write $\Lambda(x)$ for the $\cO$-lattice $\cO x + \Lambda$; 
since $\xi$ has order $N$, $[\Lambda(x):\Lambda] \ge N$.   
More generally, for any integer $m$,  
put $\Lambda(mx) = \cO \cdot m x + \Lambda = m\cO x + \Lambda$.  
Note that 
\begin{equation*}
\Lambda(mx)/\Lambda \ \cong \ \prod_{p|N} \Lambda(m x)_p/\Lambda_p 
\ = \ \prod_{p|N} (m \cO_p x + \Lambda_p)/\Lambda_p \ .
\end{equation*} 

If $p|Q$, then $G_Q$ acts on $\xi_p$ 
through the subgroup $1 + p^{\ord_p(Q)} \cO_p \subset \cO_p^{\times}$.
Noting that $\ord_p(Q_N) = \min(\ord_p(Q),\ord_p(N))$ and that
$p^{\ord_p(Q)} x \in \Lambda_p$ if $\ord_p(Q) \ge \ord_p(N)$, we have  
\begin{equation*}
G_Q \cdot \xi_p 
       \ \cong \  (x +  p^{\ord_p(Q)}\cO_p x + \Lambda_p)/\Lambda_p
       \ = \  (x +  \Lambda(p^{\ord_p(Q_N)}x)_p)/\Lambda_p \ .
\end{equation*}
Thus $\#(G_Q \cdot \xi_p) = [\Lambda(p^{\ord_p(Q_N)} x)_p:\Lambda_p]$.

If $p \!\not \, \mid Q$, then $\cO_p = \cO_{K,p}$  
and $G_Q$ acts on $\xi_p$
through  $\cO_p^{\times} \cong \prod_{\fp|p} \cO_{K,\fp}^{\times}$.  
For each $\fp|p$, and each $\cO$-lattice $L$, 
we have $L_p \cong (\cO_K L)_p$ where $\cO_K L$ is an $\cO_K$-fractional
ideal.  Thus  $\ord_{\fp}(L) := \ord_{\fp}(\cO_K L)$ is well defined.  
Write $\ord_{\fp}(\xi) = \ord_{\fp}(\Lambda) - \ord_{\fp}(\Lambda(x))$.  
Then 
$\Lambda(x)_p/\Lambda_p \ \cong \  
\prod_{\fp|p} \cO_K/\fp^{\ord_{\fp}(\xi)}$
and 
\begin{equation*}
\#(G_Q \cdot \xi_p) \ = \ [\Lambda(x)_p:\Lambda_p] 
\cdot \prod_{\substack{ \fp|p \\ 
      \ord_{\fp}(\xi) > 0 }} 
\big(1 - \frac{1}{N\fp} \big)
\end{equation*}              
where $N\fp = \#(\cO_K/\fp)$ is the norm of $\fp$.  

Combining these formulas, and using that   
$\prod_{p|N} [\Lambda(p^{\ord_p(Q_N)}x)_p : \Lambda_p] 
= [\Lambda(Q_N x):\Lambda]$, 
we obtain
\begin{equation}
\#(G_Q \cdot \xi) \ = \ [\Lambda(Q_N x):\Lambda] 
 \cdot \prod_{\substack{ \fp|N, \fp \!\not \, \mid Q, \\ 
              \ord_{\fp}(\xi) > 0 }} 
                        \big(1 - \frac{1}{N\fp} \big) \ . \label{FDV2} 
\end{equation} 

If $L$ is any $\cO$-lattice, and $F(L)$ is the area of a fundamental domain
for $\CC/L$, then by Minkowski's theorem there is a point 
$0 \ne \ell \in L$ with $|\ell| \le (4/\pi)^{1/2} F(L)^{1/2}$.  
Here, $L$ is a proper $\cO^{\prime}$-lattice 
for some order $\cO^{\prime}$ with conductor $c^{\prime}|c$.  
There are only finitely many such orders $\cO^{\prime}$, 
and for each $\cO^{\prime}$ there are only finitely many
homothety classes of proper $\cO^{\prime}$-lattices, so there are only 
a finitely many homothety classes of $\cO$-lattices.   
Hence there is a constant $C_1$, independent of $L$, such that  
$L$ has a fundamental domain $\cF(L)$ contained in the ball
$B(0,C_1 \cdot F(L)^{1/2})$.  
In turn, there is a constant $C$, independent of $L$,
such that $\cF(L) \subset \cS(0,C \cdot F(L)^{1/2})$. 
This fact is the crux of the proof.   

Again, if $L$ is an $\cO$-lattice, 
then for each ideal $\varpi$ of $\cO_K$ coprime to $c$, 
there is a unique lattice $\varpi L$  
defined by the property that $(\varpi L)_q = (\varpi \cO_K L)_q$ 
for all primes $q|N\varpi$, 
and $(\varpi L)_q = L_q$ for all primes $q \!\not \, \mid N\varpi$.  
This lattice has index $[L:\varpi L] = N\varpi$.  

\vskip .05 in
Now consider a set $\cS(a,r)$, where $a \in \CC$ and $r \le r_0$. 
For each $\sigma_i(\xi)$, we will compute 
$\#((G_Q \cdot \sigma_i(\xi))\cap \cS_E(a,r))$.  Fix $\sigma_i$, 
and replace $\xi$ by $\sigma_i(\xi)$ in the discussion above.  
Let $x \in K/\Lambda$ correspond to $\sigma_i(\xi)$, and   
let $\fp_1, \ldots, \fp_R$ be the primes of $\cO_K$ dividing $N$ but 
not $Q$, for which $\ord_{\fp}(\Lambda(x)) \ne \ord_{\fp}(\Lambda)$.   
(Note that the $\fp_j$ are independent of $\sigma_i$, since $K \subset k$
and for $p \!\not \, \mid Q$, 
$\sigma_i$ acts on $\xi$ through $\cO_p^{\times}$.)     
Then there is a one-to-one correspondence
between elements of $G_Q \cdot \sigma_i(\xi)$, and cosets $y + \Lambda$
for $y \in K$ such that $y \in x + \Lambda(Q_N x)$ 
and $y \notin \fp_j \Lambda(x)$ for $j = 1, \ldots, R$.  
Since $\Lambda(Q_Nx) \subset \Lambda(x)$, 
such $y$ necessarily belong to $\Lambda(x)$.

The lattices $\Lambda(Q_N x)$ and $\fp_1 \cdots \fp_R \Lambda(x)$ have
coprime indices in $\Lambda(x)$, so there is an $x_0 \in \Lambda(x)$ such
that 
$x_0 \equiv x \text{\ (mod $\Lambda(Q_N x)$)}$ and 
$x_0 \equiv 0 \text{ \ (mod $\fp_1 \cdots \fp_R \Lambda( x)$)}$.
Further, for any $\cO_K$-ideal $\varpi$ dividing $\fp_1 \cdots \fp_R$, 
\begin{equation*} 
\Lambda(Q_N x) \bigcap \, \big(\bigcap_{\fp_j|\varpi} \fp_j \Lambda(x) \big) 
\ = \ \varpi \Lambda(Q_N x) \ .
\end{equation*}
Clearly $y \in x+ \Lambda(Q_N x)$ if and only if $y \in x_0 + \Lambda(Q_N x)$,
and $y \in \fp_j \Lambda(x)$ if and only if 
$y \in x_0 + \fp_j \Lambda(x)$.
Since $\cS(a,r)$ injects into $\CC/\Lambda$, by inclusion-exclusion  
\begin{eqnarray}
& & \#((G_Q \cdot \sigma_i(\xi)) \bigcap \cS_E(a,r)) \label{FBB2} \\
& & \qquad \ = \ 
\sum_{\varpi|\fp_1 \cdots \fp_R} (-1)^{\lambda_K(\varpi)} \cdot 
\#( ( x_0 + \varpi \Lambda(Q_N x) ) \bigcap \cS(a,r) ) \notag
\end{eqnarray}
where $\lambda_K(\varpi)$ is the number of prime ideals of $\cO_K$ 
dividing $\varpi$.  

Take $L = \varpi \Lambda(Q_N x)$, and note that its fundamental domain 
$\cF(L)$ has area $F(L) = F \cdot N\varpi/[\Lambda(Q_N x):\Lambda]$,  
where $F$ is the area of a fundamental domain $\cF$ for $\Lambda$.  
By the same argument leading to (\ref{FDV2}) we find 
\begin{eqnarray}
& & \big|\#((x_0 + \varpi \Lambda(Q_N x)) \cap \cS(a,r)) 
 - \frac{\area(\cS(a,r))}{\area(\cF)} \cdot 
 \frac{[\Lambda(Q_N x):\Lambda]}{N\varpi} \big| \notag \\
& & \qquad \qquad \qquad \qquad  
       \ \le \  \frac{2  C S r}{F} \cdot 
       \big(\frac{[\Lambda(Q_N x):\Lambda]}{N\varpi}\big)^{1/2} 
                + \frac{ C^2S }{F} \ . \label{FCMA3}
\end{eqnarray} 

Here the index $[\Lambda(Q_N x):\Lambda]$ is independent of $\sigma_i$ by 
(\ref{FDV2}), since $\#(G_Q \cdot \sigma_i(\xi))$ and the $\fp_j$ 
are independent of $\sigma_i$.  Inserting (\ref{FCMA3}) in the 
inclusion-exclusion formula (\ref{FBB2}) and summing over all $\sigma_i(\xi)$, 
\begin{eqnarray*}
N(\xi,\cS_E(a,r)) 
& = & \frac{\area(\cS(a,r))}{\area(\cF)} \cdot 
      T [\Lambda(Q_N x):\Lambda] 
  \prod_{j=1}^R \big(1 - \frac{1}{N\fp_j} \big)  \\
& & \quad \ + \ \theta\big(\frac{2 C S r}{F} \cdot 
     T [\Lambda(Q_N x):\Lambda]^{1/2}
     \prod_{j=1}^R \big(1 + \frac{1}{N\fp_j^{1/2}} \big) \big) \\
& & \qquad    + \ \theta\big( \frac{C^2 S}{F} \cdot T 2^R \big) \ .
\end{eqnarray*}
By (\ref{FDV2}), 
$[k(\xi):k] \ = \ T [\Lambda(Q_N x):\Lambda]  \prod_{j=1}^R (1 -1/N\fp_j)$.  
Since $r \le r_0$ and $\prod_{j=1}^R(1+1/N\fp^{1/2}) \le 2^R$, 
\begin{eqnarray}
\frac{N(\xi,\cS_E(a,r))}{[k(\xi):k]} 
& = & \frac{\area(\cS(a,r))}{\area(\cF)} \notag \\
& & \quad 
 \ + \ \theta\big(\frac{2 C S r_0}{F} \cdot 
   \frac{T^{1/2} 2^R} 
            {(\prod_{j=1}^R (1 -1/N\fp_j))^{1/2}}
    \cdot \frac{1}{[k(\xi):k]^{1/2}}\big) \notag \\
& & \qquad \qquad + \ \theta\big(  \frac{ C^2 S}{F} \cdot \frac{T 2^R}
     {[k(\xi):k]} 
      \big) \ .  \label{FCN1}
\end{eqnarray}

As before $\area(\cS(a,r))/\area(\cF) = \mu(\cS_E(a,r))$.   
Here $T \le [\Gal(\kbar/k):G_Q]$ is fixed.
For each $\varepsilon > 0$ and each sufficiently large $N$,
$2^R \le 2^{\Lambda_K(N)} \le 2^{2\lambda(N)} \le d(N)^2 \le N^{\varepsilon}$.
Likewise, 
$\prod_{j=1}^R(1-1/N\fp) \ge \prod_{p|N} (1-1/p)^2 \ge C/(\log\log(N))^2$ 
for some constant $C > 0$, where the last inequality follows from 
(\cite{H-W}, Theorem 328, p.267). Finally, 
since $\xi$ has order $N$ and $Q_N \le Q$ is bounded,
$[\Lambda(Q_N x) : \Lambda] \ge N/Q$, and so 
\begin{equation} 
[k(\xi):k] \ \ge \ T \cdot N/Q \cdot C/(\log\log(N))^2 \ \ge \ 
TC/Q \cdot N^{1-\varepsilon} \label{FBMB} 
\end{equation}
for all large $N$. Combining these shows that for each $0 < \gamma < 1/2$, 
the first error term is $\cO_{\gamma}([k(\xi):k]^{-\gamma})$.
The same estimates show the second error term is 
negligible in comparison to the first.     

This completes the proof of Proposition \ref{EquiProp} 
when $E$ has complex multiplication. 

\vskip .1 in
Before closing, we note for purposes of reference 
that the arguments above provide lower bounds for the 
degree $[k(\xi):k]$ in terms of the order $N$ of $\xi$.  When $E$ does not have
complex multiplication, then since $T$ is fixed, $Q_N \le Q$, 
and $\prod_p (1-1/p^2)$ converges to a nonzero limit, (\ref{FGJ1}) shows there
is a constant $C_1$ depending only on $E$ such that 
\begin{equation} \label{FLL1} 
[k(\xi):k] \ \ge \ C_1 N^2  \ .
\end{equation}
When $E$ has complex multiplication, then since $T$ and $Q$ are fixed, (\ref{FBMB})
shows that there is a constant $C_2$ depending only on $E$ such that 
\begin{equation} \label{FLL2} 
[k(\xi):k] \ \ge \ C_2 N/(\log\log(N))^2 \ .
\end{equation} 
 
\end{proof}

\vskip .3 in


\vskip .2 in

\begin{thebibliography}{DHK02}

\bibitem[Ba75]{Baker} A. Baker, Transcendental number theory, 
Cambridge University Press, Cambridge, 1975.

\bibitem[BHpp]{B-H} M. Baker and L. C. Hsia, Canonical heights, transfinite diameters, 
and polynomial dynamics, to appear in J. Reine Angew. Math.

\bibitem[B1886]{Bang} A. Bang, Taltheoreske Unders{\o}gelser, 
Tidsskr. Math. (5) 4 (1886), 70-80 and 130-137.

\bibitem[CH99]{Ch-Han} J. Cheon and S. Hahn, The orders of the reductions of a point 
in the Mordell-Weil group of an elliptic curve, Acta Arithmetica 88 (1999), 219-222.

\bibitem[CLpp]{CL} A. Chambert-Loir, Mesures et \'equidistribution sur les
espaces de Berkovich, preprint, arXiv:math.NT/0304023 v3.  

\bibitem[Co86]{Conway} J.  Conway, Functions of one complex variable 
($2^{nd}$ edition), Springer-Verlag, New York, 1986.

\bibitem[DHK02]{D-HK} S. David and N. Hirata-Kohno, Recent progress on linear forms
in elliptic logarithms, pp.26-37 in:  
G. W\"ustholz, ed., A panorama of number theory, 
or the view from Baker's garden, Cambridge University Press, Cambridge, 2002.

\bibitem[EW99]{E-W} G. Everest and T. Ward, Heights of polynomials and entropy in 
algebraic dynamics, Springer-Verlag, New York, 1999.  

\bibitem[EF96]{EnF} G. Everest and B. N\'i Fhlath\'uin, The elliptic Mahler measure,
Math. Proc. Cambridge Philos. Soc. 120 (1996), 13-25.

\bibitem[FRLpp]{F-RL} C. Favre and J. Rivera-Letelier, Equidistribution quantitative
des points de petite hauteur sur la droite projective, preprint, 2005.  

\bibitem[HW71]{H-W} G. Hardy and E. Wright, An introduction to the theory of numbers 
($4^{th}$ edition), Oxford University Press, London, 1971.

\bibitem[L73]{Lang1} S. Lang, Elliptic functions, Addison-Wesley, Reading, 1973.

\bibitem[Ra83]{Ray} M. Raynaud, 
Courbes sur une vari\'et\'e ab\'elienne et points de torsion,
Inv. Math 71 (1983), 207-233. 

\bibitem[Ru89]{R1} R. Rumely, Capacity theory on algebraic curves, Lecture Notes
in Mathematics 1378, Springer-Verlag, Berlin-Heidelberg-New York, 1989.
 
\bibitem[Sc74]{Schinzel} A. Schinzel, Primitive divisors of the expression $A^n-B^n$
in number fields, J. Reine Angew. Math 268 (1974), 27-33.  

\bibitem[Se72]{Serre} J.-P. Serre, Propri\'et\'es galoisiennes des points d'order 
fini des courbes elliptiques, Invent. Math 15 (1972), 259-331.

\bibitem[Sh71]{Sh} G. Shimura, Introduction to the arithmetic theory 
of automorphic functions, Princeton University Press, 
USA, 1971.

\bibitem[Si86]{Sil1} J. Silverman, The arithmetic of elliptic curves, 
Graduate Texts in Mathematics 106, Springer-Verlag,
Berlin-Heidelberg-New York, 1986. 

\bibitem[Si99]{Sil2} J. Silverman, Advanced topics in the theory of elliptic curves,
Graduate Texts in Mathematics 151, Springer-Verlag,
Berlin-Heidelberg-New York, 1999. 

\bibitem[Si95]{Sil3} J. Silverman, 
Exceptional units and numbers of small Mahler measure,
Experimental Mathematics 4 (1995), 69-83.  

\bibitem[STpp]{SZ-T} L. Szpiro and T. Tucker, 
Equidistribution and generalized Mahler measures, preprint, 2005.  

\bibitem[U95]{Ull0} E. Ullmo, 
Points entiers, points de torsion et amplitude arithm\'etique,
American Journal of Mathematics 117 (1995), 1039-1055.   

\bibitem[U98]{Ull} E. Ullmo, Positivit\'e et discr\`etion des points alg\'ebriques 
des courbes, Annals of Mathematics (2) 147 (1998), 167-179.  

\end{thebibliography}
\end{document}